\documentclass[12pt,a4paper,final]{iopart}

\usepackage{iopams}  
\usepackage[font=footnotesize,labelfont=bf,margin=2cm]{caption}
\usepackage{cite}

\expandafter\let\csname equation*\endcsname\relax
\expandafter\let\csname endequation*\endcsname\relax
\newtheorem{theorem}{Theorem}

\newtheorem{prop}[theorem]{Proposition}
\newtheorem{lemma}[theorem]{Lemma}
\newtheorem{remark}[theorem]{Remark}
\usepackage{amsmath}
\usepackage{amssymb}
\usepackage{mathtools}
\usepackage{verbatim}
\newcommand{\floor}[1]{\lfloor #1 \rfloor}
\usepackage{graphicx}
\usepackage[breaklinks=true,colorlinks=true,linkcolor=blue,urlcolor=blue,citecolor=blue]{hyperref}
\def\E{{\mathbb E}}

\def\R{{\mathbb R}}

\begin{document}
\title[Preparing an article for IOP journals in  \LaTeXe]{Exact and asymptotic properties of $\delta$-records in the linear drift model.}

\author{R. Gouet$^{1}$, M. Lafuente$^{2}$, F. J. L\'opez$^{2,3}$ and G. Sanz$^{2,3}$}
\address{$^1$Departamento de Ingenier\'ia Matem\'atica and CMM (UMI 2807, CNRS), Universidad de Chile, Avenida Blanco Encalada 2120,837-0459, Santiago, Chile.}
\address{$^2$Departamento de M\'etodos Estad\'isticos, Facultad de Ciencias, Universidad de Zaragoza, C/ Pedro Cerbuna, 12, 50009 Zaragoza, Spain.}
\address{$^3$Instituto de Biocomputaci\'on y F\'isica de Sistemas Complejos (BIFI), Universidad de Zaragoza, 50018 Zaragoza, Spain.}

\eads{\mailto{rgouet@dim.uchile.cl}, \mailto{miguellb@unizar.es}, \mailto{javier.lopez@unizar.es}, \mailto{gerardo.sanz@unizar.es}}
\begin{abstract}
The study of records in the Linear Drift Model (LDM) has attracted much attention recently due to
applications in several fields. In the present paper we study
$\delta$-records in the LDM, defined as observations which
are greater than all previous observations, plus a fixed real quantity
$\delta$. We give analytical properties of the probability of
$\delta$-records and study the correlation between $\delta$-record
events. We also analyse the asymptotic behaviour of the number of
$\delta$-records among the first $n$ observations and give conditions
for convergence to the Gaussian distribution. As a consequence of
our results, we solve a conjecture posed in J. Stat. Mech. 2010, P10013, regarding the total
number of records in a LDM with negative drift. Examples of
application to particular distributions, such as Gumbel or Pareto are
also provided. We illustrate our results with a real data set of
summer temperatures in Spain, where the LDM is consistent with the
global-warming phenomenon.

\end{abstract}

\vspace{2pc}

\tableofcontents

\markboth{Exact and asymptotic properties of $\delta$-records in the linear drift model.}{}

\section{Introduction.}

Extreme values  and records  have attracted large efforts and attention since the beginnings of statistics and probability, due to their intrinsic interest and their mathematical   challenges. An important motivation for studying records comes from their connections with other interesting problems and, of course, from their countless practical applications in different fields such as climatology \cite{Benestad,Majumdar,Redner,Wergen2}, sports \cite{Einmahl,Gembris1,Gembris2}, finance \cite{Wergen4,Wergen0} or biology \cite{Park}. Moreover, records have been used in statistical inference because, in some contexts, data is inherently composed of record observations \cite{Dey,Gulati,Jafari,Wang}. The classical probabilistic setting of independent and identically distributed random observations (iid)  observations has been profusely studied. Main results in this framework can be found in the  monographs  \cite{Ahsanullah,Arnold,Nevzorov}. In the last years there has been an increasing interest in the study of records in correlated observations such as random walks or time series \cite{Godreche_moving,Godreche1,Godreche2,GouetAAP,Kearney_2020,Mounaix}. 

An interesting departure from the iid model, which introduces time-dependence between observations, results from adding a deterministic linear trend to the iid observations, thus obtaining the so named Linear Drift Model (LDM). This model was first introduced in \cite{Ballerini-Resnick-1} and later developed in \cite{Ballerini-Resnick-2,Borovkov1,DeHaan}. The model was also considered in \cite{Franke2}, under a wide range of scenarios, and has proven particularly useful in the study of global warming phenomena \cite{Wergen2,Wergen3}. Furthermore, the importance of this model is not only related to applications but also to its mathematical structure. For instance, the study of records in the LDM model can be helpful in determining whether the underlying distribution is heavy-tailed or not \cite{Franke1,Wergen1}.

Also, different generalizations of the notion of record, such as  
near-records \cite{Balakrishnan,Hashorva,Pakes} or $\delta$-exceedance records \cite{Balakrishnan_exceedance,Park_exceedance} have been proposed recently. We will work with $\delta$-records, first introduced in \cite{Gouet1}, which are observations greater than all previous entries, plus a fixed quantity $\delta$. In the iid setting, the distribution \cite{LopezSalamanca1,LopezSalamanca2}, process structure \cite{Gouet3} and asymptotic properties \cite{Gouet2} of $\delta$-records have been studied. In the case $\delta<0$, where $\delta$-records are more numerous than records, their  use in statistical inference has been recently proposed and positively assessed; see \cite {Gouet2,Gouet0,Gouet222}. 

In this work, we study  $\delta$-records from observations obeying the LDM, while revisiting some open questions about records. We  analyse the positivity and continuity of the asymptotic $\delta$-record probability as a function of $\delta$ and of the trend parameter $c$. We also obtain a law of large numbers and a central limit theorem for the counting process of $\delta$-records, thus extending the corresponding results  in \cite{Ballerini-Resnick-1}. Furthermore, we completely characterize the finiteness of the number of $\delta$-records and, in particular, we solve a conjecture posed in \cite{Franke2}, about the finiteness of the number of usual records in the LDM with negative trend. 

We assess the effect of $\delta$ on the $\delta$-record probabilities and correlations, for explicitly solvable models. Some of the results obtained in these examples are new and shed light on the behaviour of record events, when the underlying distribution is heavy-tailed. Finally we illustrate our results by analyzing a real dataset of temperatures, which fits the LDM with a trend parameter consistent with the global-warming phenomenon.

\section{$\delta$-records in the linear drift model}

Our objects of interest in this paper are $\delta$-records, formally defined as follows: given a sequence of observations $(Y_n)_{n\ge1}$ and $\delta\in\mathbb{R}$ a parameter, $Y_1$ is defined conventionally as $\delta$-record and, for $j\ge2$, $Y_j$ is a $\delta$-record if $Y_j>\max\{Y_1,\ldots,Y_{j-1}\} +\delta$.

Note that  $\delta$-records are just (upper) records, if $\delta=0$. If $\delta>0$, a $\delta$-record is necessarily a record and  $\delta$-records are a subsequence of records. On the other hand, if $\delta<0$, a $\delta$-record can be smaller than the current maximum, so records are a subsequence of $\delta$-records.

Throughout this paper we assume that the $Y_n$ are random variables obeying the LDM, that is, $Y_n$ can be represented as
\begin{equation}
\label{lineardrift}
Y_n=X_n+cn,\,n\ge1,
\end{equation}
where $c\in\mathbb{R}$ is the trend parameter and $(X_n)_{n\ge1}$ is a sequence of iid random variables, with (absolutely continuous) cumulative distribution function (cdf) $F$ and probability density function (pdf) $f$.  Another important parameter of the model is the right-tail expectation of the $X_j$, defined as
\begin{equation*}
\mu^+=\int_{0}^{\infty}xf(x)dx.
\end{equation*}
For simplicity, we assume the existence of an interval of real numbers $I=(x_-,x_+)$, with $-\infty\le x_-<x_+\le\infty$, such that $f(x)>0$, for all $x\in I$, and $f(x)=0$ otherwise. Note that $x_-=\inf\{x:F(x)>0\}$ and $x_+=\sup\{x:F(x)<1\}$.

Let $1_{j,\delta}$ denote the indicator of the event $\{Y_j$ is a $\delta$-record$\}$. That is,  $1_{j,\delta}=1$ if $Y_j>\max\{Y_1,\ldots,Y_{j-1}\} +\delta$ and $1_{j,\delta}=0$ otherwise. So, the number of $\delta$-records up to index $n$ is computed as $N_{n,\delta}=\sum_{j=1}^n 1_{j,\delta}$. 

Under the LDM, the probability of \{$Y_j$ is a $\delta$-record\} is easily computed by conditioning, as 
\begin{equation*}
p_{j,\delta}:=\mathbb{E}[1_{j,\delta}]=\int_{-\infty}^\infty \prod_{i=1}^{j-1}F(x+ci-\delta)f(x)dx,\label{eq:pndelta}
\end{equation*}
where $\mathbb{E}[\cdot]$ denotes the mathematical expectation.
Moreover, the asymptotic $\delta$-record probability is given by the formula
\begin{equation}
\label{eq:pdelta}
p_{\delta}:=\lim_{n\to \infty}p_{n,\delta}=\int_{-\infty}^\infty \prod_{i=1}^{\infty}F(x+ci-\delta)f(x)dx,
\end{equation}
which is mathematically justified by the monotone convergence theorem for integrals. 

In what follows we occasionally write $1_{j,\delta}(c), N_{n,\delta}(c), p_{j,\delta}(c), p_\delta(c)$, etc.  to emphasize the dependence on the trend parameter $c$.

\section{Properties of the $\delta$-record probabilities}\label{seccionanalitica}

We begin with a simple property about the asymptotic $\delta$-record probability of an affine transformation of the LDM. 
Let $\tilde{X}_n=bX_n+a$, with $b>0$, $a\in\mathbb{R}$, and  $\tilde{Y}_n=\tilde{X}_n+cn, n\ge 1$. If  $\tilde{p}_\delta(c)$ is  the $\delta$-record probability in this model, then it holds 
\begin{equation*}
\tilde{p}_\delta(c) =p_{\frac{\delta}{b}}(\tfrac{c}{b}).
\end{equation*}

We consider next some analytical properties of $p_{j,\delta}(c)$ and $p_{\delta}(c)$, as functions of $c$ and $\delta$. We note first that both are increasing in $c$ and decreasing in $\delta$. Moreover, it is easy to see that $p_{j,\delta}(c)$ is decreasing in $j$ and continuous in $c$, converging  to 1 as $c\to\infty$. The continuity of $p_{\delta}(c)$ is less clear because of the infinite product within the integral in \eqref{eq:pdelta}.
\subsection{Positivity of $p_\delta(c)$}\label{posp}
We show that the positivity of $p_\delta(c)$ depends on $c$ and $\delta$ and on the right-tail behaviour of $F$. We consider two cases depending on  $\mu^+$:

\noindent {\bf 1}. $\mu^+=\infty$. In this case $p_\delta(c)=0$, for all $\delta,c\in\R$. 

To justify this claim, we show that  $\prod_{j=1}^\infty F(x+cj-\delta)=0,$ for all $x\in(x_-,x_+)$. 

If $c<0$ the conclusion is immediate because $F(x+cj-\delta)\to 0$, as $j\to\infty$. 

If $c=0$, we note that $\mu^+=\infty$ implies $x_+=\infty$ and so,  $F(x-\delta)<1$. Thus $\prod_{j=1}^\infty F(x+cj-\delta)=0$. 

Finally, if $c>0$, we note that $\mu^+=\infty$ implies $\sum_{i=1}^{\infty}(1-F(x+ci-\delta))=\infty$,  which in turn implies  $\prod_{j=1}^\infty F(x+cj-\delta)=0$. This follows from the definition of $\mu^+$ and from Taylor's expansion of $\log(1+x)$.

Distributions with $\mu^+=\infty$ can be considered as ``right-heavy-tailed'' and we observe  that, for such distributions, the linear trend has no impact on the asymptotic probability of a $\delta$-record. This class of distributions includes the Pareto and Fr\'echet,  with shape parameter $\alpha\in(0,1]$. 

\noindent {\bf 2}. $\mu^+<\infty$. As in the previous case, we have three situations depending on the sign of $c$. 

For $c<0$,  $p_\delta(c)=0$, for all $\delta\in\mathbb{R}$,  since $\prod_{j=1}^\infty F(x+cj-\delta)= 0$, for all $x\in(x_-,x_+)$.

 If $c=0$, 
\begin{equation}\label{dosest}
p_\delta(0)=\int_{-\infty}^{\infty}\prod_{j=1}^\infty F(x-\delta)f(x)dx=\int_{x_++\delta}^\infty f(x)dx,
\end{equation}
which is positive if and only if $x_+<\infty$ and $\delta<0$.

 Finally, if $c>0$, then $p_\delta(c)=0$ if  and only if $x_+-x_-\le\delta-c$. Indeed, note that, if $x_+-x_-\le\delta-c$, then $\mathbb{P}[Y_n>Y_{n-1}+\delta]=0$, for all $n$,  and so, only the first observation (by convention) is a $\delta$-record. Conversely, if $x_+-x_->\delta-c$,  then the interval $J:=(x_-,x_+)\cap(x_--c+\delta,\infty)$ is nonempty and, for every $x\in J$, we have $F(x+cj-\delta)\ge F(x+c-\delta)>0$, for all $j$. Now, since $F(x+cj-\delta)\to1$ as $j\to\infty$, and $\mu^+<\infty$, we have $\sum_{j=1}^\infty(1-F(x+cj-\delta))<\infty$, which implies $\prod_{j=1}^\infty F(x+cj-\delta)>0$ and, so $p_\delta(c)>0$.

\vspace{3pt}
Summarizing the above findings, we state
\begin{theorem}\label{thm:positivity}
	$p_\delta(c)>0$ if and only if $\mu^+<\infty$ and one of the following conditions holds
	\begin{enumerate}
		\item $c>0\ and\ \delta<x_+ - x_- +c$,
		\item $c=0,\ \delta <0 \ and \ x_+<\infty$.
	\end{enumerate}
\end{theorem}
\subsection{Continuity of \texorpdfstring{$p_\delta(c)$}.}
As commented at the beginning of this section, the continuity of $p_\delta(c)$ is not obvious. However, thanks to Theorem \ref{thm:positivity} we can restrict attention to distributions $F$ with finite right-tail expectation since, otherwise, $p_\delta(c)$ vanishes and continuity is trivial. Thus, we assume throughout this section that $\mu^+<\infty$.

A first interesting fact, which is rigorously proved in  Proposition \ref{prop:continuidad1} of the Appendix, is that $\prod_{i=1}^\infty F(x+ci-\delta)$ is continuous at every $c\ne0$, for every $x\in(x_-,x_+)$, such that $x\ne x_-+\delta-c$. Then, thanks to the bounded convergence theorem of integration, we conclude that $p_\delta(c)$ is continuous, at every $c\ne0$.

The continuity at $c=0$ is subtler to establish and depends of the sign of $\delta$ and the finiteness of $x_+$, the right-end point of $F$. 
Note that, for every $c>0$ and $N\ge1$, we have
\begin{equation*} \prod_{j=1}^\infty F(x-\delta)\le\prod_{j=1}^\infty F(x+cj-\delta)\le\prod_{j=1}^N F(x+cj-\delta).\end{equation*}
Then, taking the limit as $c\to0^+$ in the above inequalities, 
\begin{equation*} \prod_{j=1}^\infty F(x-\delta)\le\lim_{c \to 0^+}\prod_{j=1}^\infty F(x+cj-\delta)\le F(x-\delta)^N. \end{equation*} Therefore, $\lim_{c \to0^+}\prod_{j=1}^\infty F(x+cj-\delta)$ is 0, if $x<x_++\delta$, and 1 otherwise. 
Then, by the dominated convergence theorem,
\begin{equation*}\label{mct}
\lim_{c\to 0^+} p_\delta(c)=
 \int_{-\infty}^{\infty}\lim_{c \to0^+}\prod_{j=1}^\infty F(x+cj-\delta)f(x)dx=\int_{x_++\delta}^\infty f(x)dx.
\end{equation*}
Thus, $p_\delta(c)$ is right-continuous at $c=0$  by \eqref{dosest}.
Regarding left-continuity at 0, recall that $p_\delta(c)=0$ for $c<0$. So, $p_\delta(c)$ is discontinuous at 0 if and only if $x_+<\infty$ and $\delta<0$. 

We now show the continuity of  $p_{\delta}(c)$ as a function of $\delta$. The result is trivial if $c<0$, since $p_{\delta}(c)=0$, for all $\delta \in \mathbb{R}$. For $c=0$, note that, by \eqref{dosest}, $p_\delta(0)=1-F(x_++\delta)$, which is continuous since $F$ is a continuous function.

If $c>0$ and $(\delta_n)_{n\ge1}$ is a sequence converging to $\delta$, we prove that
\begin{equation}\label{limitedelta}
\lim_{n\to\infty}\prod_{i=1}^{\infty}F(x+ci-\delta_n)=\prod_{i=1}^{\infty}F(x+ci-\delta),
\end{equation}
for all $x\in(x_-,x_+)$, $x\ne x_-+\delta-c$. Indeed, let $x<x_-+\delta-c$, then $F(x+c-\delta)=0$ yielding $\prod_{i=1}^{\infty}F(x+ci-\delta)=0$. Also $F(x+c-\delta_n)=0$ for $n$ large enough and \eqref{limitedelta} follows. Let now $x>x_-+\delta-c$ and $\varepsilon>0$ such that $x+c-\delta-\varepsilon>x_-$. Then, for $n$ large enough, we have $\mid\delta_n-\delta\mid <\varepsilon$ and
\begin{align*}
-\sum_{i=1}^\infty \log F(x+ci-\delta_n)\le -\sum_{i=1}^\infty \log F(x+ci-(\delta+\varepsilon))<\infty,
\end{align*}
since $\mu^+<\infty$. So \eqref{limitedelta} holds, and continuity follows.

In the following theorem we summarize conditions for continuity of $p_\delta(c)$.
\begin{theorem}\label{tmacontinuity}\ 
	The asymptotic $\delta$-record probability $p_\delta(c)$, as a function of $c,\delta$, is
	\begin{itemize}
		\item[(a)]  continuous at every  $c\ne 0$ and right-continuous at $c=0$, for all $\delta$;
		 \item[(b)] discontinuous at $c=0$ if and only if $x_+<\infty, \delta<0$, and
		\item[(c)]  continuous in $\delta$, for all $c$.
	\end{itemize}
\end{theorem}

\section{Exactly solvable models}\label{seccion_ejemplos}

In general  it is not possible to compute exactly the probabilities $p_{j,\delta}$ or $p_\delta$. We show below explicit results for the Gumbel distribution and for particular instances of the Dagum family of distributions.
\subsection{The Gumbel distribution}\label{Gumbelsection}
Let  $F(x)=\exp(-\exp(-x))$, for $x\in\R$, be the Gumbel distribution. Note that $F(x+cj-\delta)=F(x)^{e^{-cj+\delta}}$. Then, if $c\ne0$, \begin{equation*}\prod_{j=1}^{n-1} F(x+cj-\delta)=F(x)^{\sum_{j=1}^{n-1}e^{-cj+\delta}}=F(x)^{e^\delta\frac{e^{-c}-e^{-nc}}{1-e^{-c}}}\end{equation*} and, if $c=0$, $\prod_{j=1}^{n-1} F(x+cj-\delta)=F(x)^{(n-1)e^{\delta}}$.
So, from \eqref{eq:pndelta} we get
\begin{equation*}
\begin{split}
p_{n,\delta}(c)&=\int_{-\infty}^{\infty}F(x)^{e^\delta\frac{e^{-c}-e^{-nc}}{1-e^{-c}}}f(x)dx\\
&=\frac{1-e^{-c}}{1-e^{-c}+e^\delta(e^{-c}-e^{-nc})},
\end{split}
\end{equation*}
if $c\ne0$, and
\begin{equation*}
p_{n,\delta}(0)=\frac{1}{(n-1)e^{\delta}+1}.
\end{equation*}
Note that, taking limits as $n\to\infty$, in the above formulas, we obtain 
\begin{equation*}
p_{\delta}(c)=\frac{1-e^{-c}}{e^\delta e^{-c}+1-e^{-c}}=\frac{1}{1+\frac{e^{-c}}{1-e^{-c}}e^\delta},
\end{equation*}
if $c>0$ and $p_\delta(c)=0$, if $c\le0$, as expected from Theorem \ref{thm:positivity}. 

\begin{figure}
	\centering
	\includegraphics[width=0.5\textwidth]{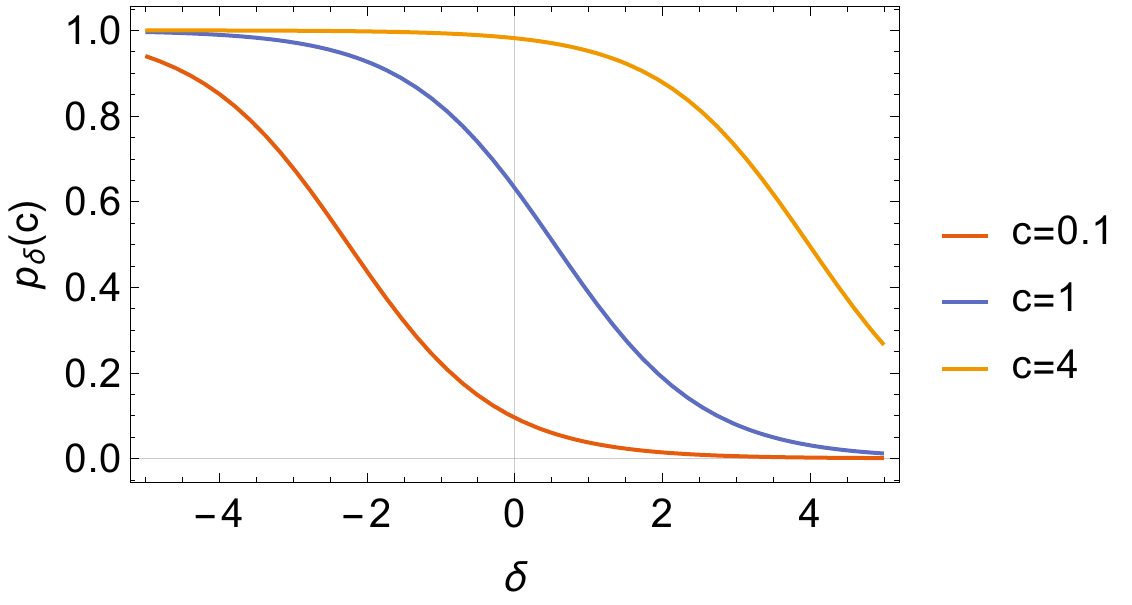}\includegraphics[width=0.5\textwidth]{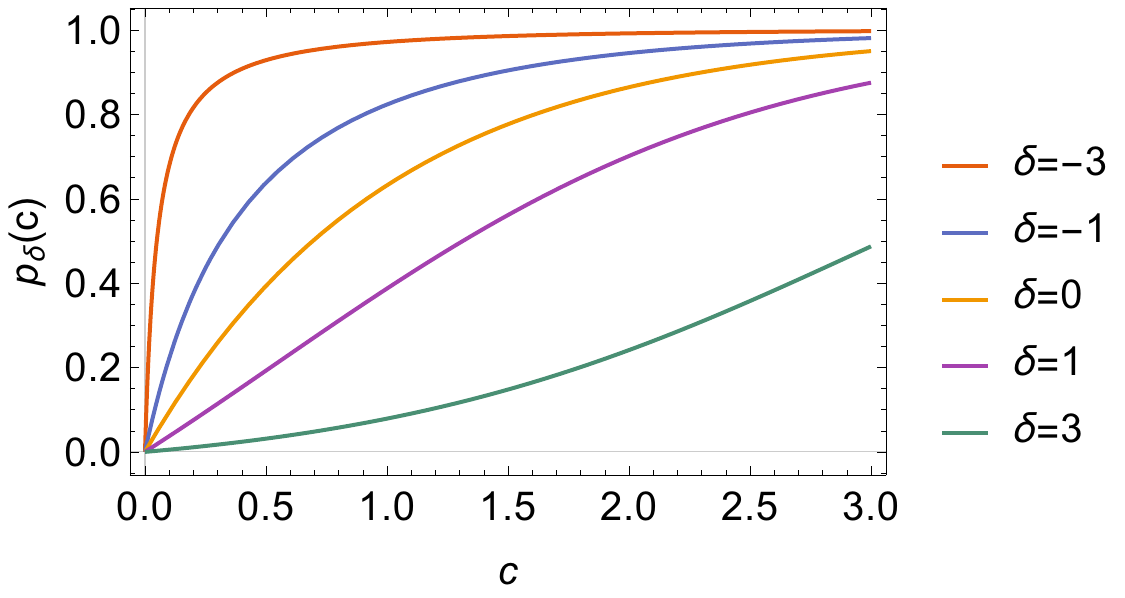}
	\caption{Asymptotic $\delta$-record probability $p_\delta(c)$ for the Gumbel distribution as function of $\delta$ and $c$.}
	\label{plotProbGumbel}
\end{figure}

Also, for every $c>0$,  $p_\delta(c)$ decreases with $\delta$ as a logistic function of $-\delta$. Figure \ref{plotProbGumbel} shows the behaviour of $p_\delta(c)$ as a function of $\delta$ and $c$. 

\subsection{The Dagum family of distributions}\label{Dagumsection} 
The random variables in the Dagum family of distributions have cdf given by $F(x)=\left(1+\left(\frac{x}{b}\right)^{-a}\right)^{-q}$, for  $x\ge0$, and $F(x)=0$, for $x<0$, where $a,b,q$ are positive parameters. Two important distributions in the Dagum family are the Loglogistic, with parameters $a,b>0$, $q=1$, and the Pareto (up to a shift), with $a>0$, $b=q=1$. For simplicity, in this example we limit our attention to the case $a=1$, which has $\mu^+=\infty$. 

By Theorem \ref{thm:positivity} we know that $p_\delta(c)=0$, for every $c,\delta\in\mathbb{R}$, so we chose to  analyse the speed of convergence of $p_{n,\delta}(c)$ to 0, for some values of $c,\delta$. To that end, observe that
the formula for $p_{n,\delta}$ takes the manageable form
\begin{equation}
\label{eq:pndeltaDagum}
p_{n,\delta}(c)=\int_{(\delta-c)^+}^\infty\prod_{i=1}^{n-1}\left(\frac{x+ci-\delta}{x+b+ci-\delta}\right)^qf(x)dx,
\end{equation}
which becomes simpler if we further assume that $c=b$ (that is, the trend parameter of the LDM is equal to the scale parameter of the distribution). From \eqref{eq:pndeltaDagum} we get
\begin{equation}
p_{n,\delta}(c)=\int_{(\delta-c)^+}^\infty\left(\frac{x+c-\delta}{x+cn-\delta}\right)^qf(x)dx.\label{eq:pndeltaDagumb=c}
\end{equation}
Note that the Pareto(1,1) distribution, taking  $c=1$, is included as a particular case. This distribution will be studied at the end of this example and later,  in section \ref{sectionCorrelations}, in the context of $\delta$-record correlations.

We introduce the notation $p_{n,\delta}^{(q)}(c)$ to make explicit the dependence of $p_{n,\delta}(c)$ on $q$.
First, for records ($\delta=0$) we have, 
\begin{align}
p_{n,0}^{(q)}(c)=&
cq\int_{0}^\infty x^{q-1}\left( x+cn \right)^{-q}(x+c)^{-1}dx \nonumber \\
=&qn^{-q}\int_0^1 t^{q-1}(1-t(n-1)/n)^{-q}dt\label{preHyper1}\\
=&\frac{q}{(n-1)^q}\int_1^n\frac{(y-1)^{q-1}}{y}dy, \label{preHyper2}
\end{align}
where the second equality follows from the change of variable $x=ct/(1-t)$ and the third from $1-t(n-1)/n=1/y$.

Observe that \eqref{preHyper1} and \eqref{preHyper2} do not depend on $c$ and so, for the sake of simplicity, we  write $p_{n,0}^{(q)}$. Moreover, from formula \eqref{preHyper1} we see that $$p_{n,0}^{(q)}=n^{-q}\ _2 F_1\left(q,q;q+1;(n-1)/n\right),$$ where $_2F_1$ is the Gauss hypergeometric function. 

Also, from \eqref{preHyper2} and  using the binomial expansion, for $q=1,2,\ldots$, we readily obtain
\begin{equation}\label{dagumexplicita}
p_{n,0}^{(q)}=\frac{q}{(n-1)^q}\left((-1)^{q-1}\log n+
\sum_{k=1}^{q-1}{q-1\choose k}\frac{(-1)^{q-1-k}}{k}(n^k-1)\right).
\end{equation}

The asymptotic behaviour of $p_{n,0}^{(q)}$, for any $q\in(0,\infty)$, can be obtained from \eqref{preHyper2}.  For $q=1$, \eqref{dagumexplicita} yields $p_{n,0}^{(1)}=\frac{1}{n-1}\log n$. 
For $q>1$, the leading term in the integral in \eqref{preHyper2} is $y^{q-2}$, so $p_{n,0}^{(q)}\sim \frac{q}{q-1}\frac{1}{n}$. For $q\in(0,1)$, the integral in \eqref{preHyper2} converges and, using formula 3.191.2 in \cite{GR}, we get
\begin{equation*}
p_{n,0}^{(q)}\sim n^{-q}q\int_1^\infty\frac{(y-1)^{q-1}}{y}dy=n^{-q}q\Gamma(1-q)\Gamma(q).
\end{equation*}
Thus,
\begin{align}
p_{n,0}^{(q)} \sim \begin{cases} n^{-q}q\Gamma(1-q)\Gamma(q), & \mbox{if } 0<q<1,\\
\log(n)/n, & \mbox{if } q=1,\\ n^{-1}\frac{q}{q-1}, & \mbox{if } q>1.\\ \end{cases}\label{pnbeha}
\end{align}

It is interesting to observe that the limiting behaviour of $p_{n,0}^{(q)}$, as a function of the power of the tail $q$, seems to match the asymptotic behaviour of $p_{n,0}(c)$ when $F$ is the  Fr\'echet distribution ($F(x)=\exp(-x^{-1})$, $x>0$) and the tuning parameter is the trend $c$, studied in \cite{DeHaan}.

We now consider $\delta\ne0$ and investigate whether $p_{n,\delta}^{(q)}/p_{n,0}^{(q)}\to1$, as $n\to\infty$. This result can be expected since, as $\mu^+=\infty$, the variables $X_n$ take very large values, so $\delta$ may have little influence on the probability of $\delta$-record, in the long term.

From (\ref{eq:pndeltaDagumb=c}) we may evaluate $p_{n,\delta}^{(q)}$, for any $q\in\mathbb{N}$, although the computation becomes lengthy as $q$ grows. We have carried out the computation with values of $q$ from 1 to 7, and obtained
\begin{equation*}
p_{n,\delta}^{(1)}\sim \frac{\log(n)}{n},\qquad
p_{n,\delta}^{(q)}\sim \frac{q}{q-1}\frac{1}{n},\quad q=2,\ldots,7.
\end{equation*}
So, from \eqref{pnbeha} we have  $p_{n,\delta}^{(q)}/p_{n,0}^{(q)}\to1$, at least for $q=1,\ldots,7$.

For noninteger values of $q\in(0,\infty)$, the limit behaviour of (\ref{eq:pndeltaDagumb=c}) is harder to analyse. To get a tractable expression, we impose $\delta=c$. Proceeding as above, we have, for $n>2$,
\begin{equation*}\label{deltaigualc}
p_{n,\delta}^{(q)}=\frac{q(n-1)^q}{(n-2)^{2q}}\int_{1}^{n-1}\frac{(y-1)^{2q-1}}{y^{q+1}}dy.\end{equation*}
Therefore, we have
\begin{align*}
p_{n,\delta}^{(q)} \sim \begin{cases} n^{-q}\frac{\Gamma(2q)\Gamma(1-q)}{\Gamma(q)}, & \mbox{if } 0<q<1,\\
\log(n)/n, & \mbox{if } q=1,\\ n^{-1}\frac{q}{q-1}, & \mbox{if } q>1.\\ \end{cases}
\end{align*}
So, under the above stated conditions, $p_{n,\delta}^{(q)}\sim p_{n,0}$, for $q\ge1$, but this is not the case if $q\in(0,1)$. 

\begin{figure}
	\centering
	
	  \includegraphics[width=0.5\textwidth]{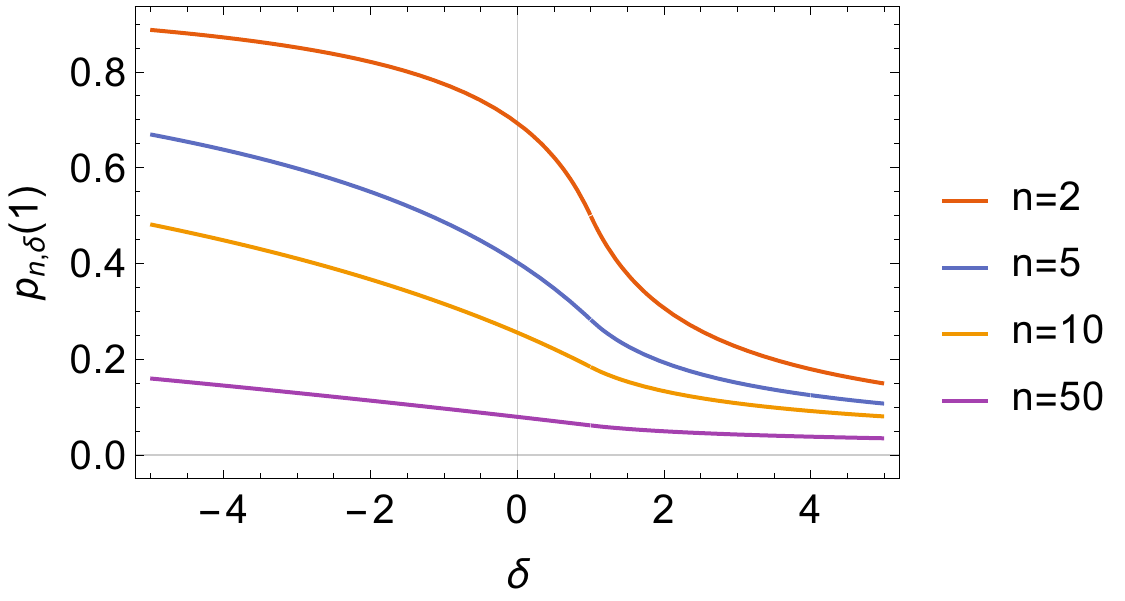}\includegraphics[width=0.5\textwidth]{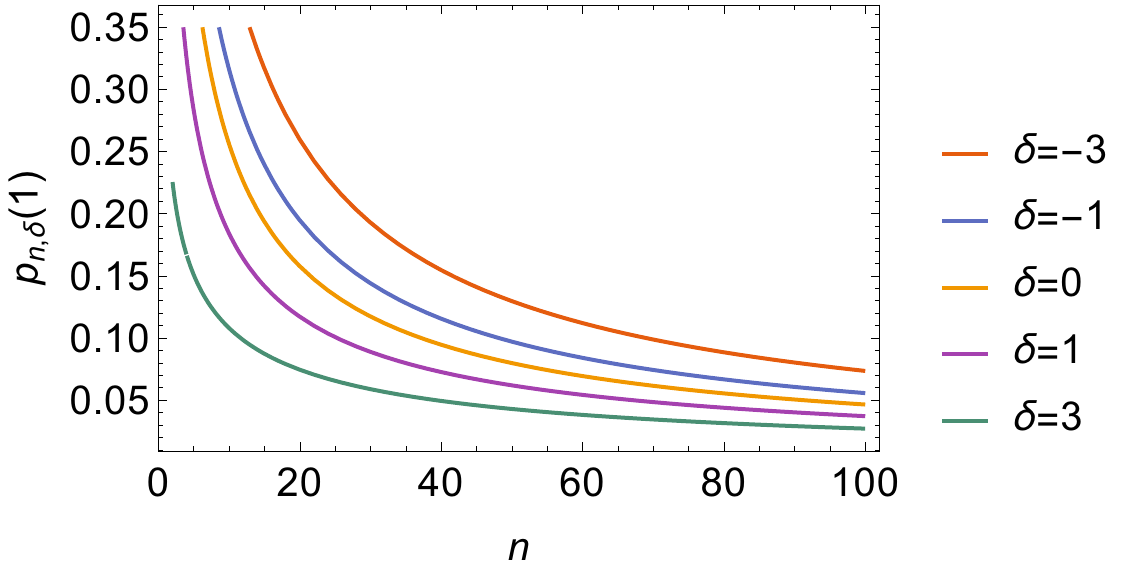}
	\caption{$\delta$-record probability $p_{n,\delta}(c)$ for the Pareto Distribution as a function of $\delta$ and $n$ with $c=1$.}
	\label{plotProbPareto}
\end{figure}

To conclude this example, we study the particular case of the Pareto distribution, that is, $F(x)=\left( 1-1/x\right)1_{\{x>1\}}$, and take $c=1$. The probability of $\delta$-record is explicitly computed as:
\begin{equation}\label{probpareto}\begin{split}p_{n,\delta}&=\int_{\delta\vee1}^\infty
\frac{x-\delta}{x^2(x+n-1-\delta)}dx\\
&=\frac{1}{(n-1-\delta)^2}\left((n-1)\log(\tfrac{n-\min\{1,\delta\}}{\max\{1,\delta\}})-\min\{1,\delta\}(n-1-\delta)\right),\end{split}\end{equation}
if $\delta\ne n-1$ and $p_{n,\delta}=\frac{1}{2(n-1)}$, if $\delta=n-1$. Figure \ref{plotProbPareto} shows the behaviour of $p_{n,\delta}$ as a function of $n$ and $\delta$.

\section{Correlations}\label{sectionCorrelations}

The indicators of $\delta$-records are in general not independent in the case of  iid random variables, see \cite{Gouet2}. In \cite{Wergen1} the authors study the dependence of record events in the LDM, by means of the following dependence index ($\delta=0$ in their case)
\begin{equation}
l_{n}(c,\delta):=\frac{\mathbb{P}[\text{obs. }n\text{ and }n+1\text{ are }\delta\text{-records}]}{\mathbb{P}[\text{obs. }n\text{ is }\delta\text{-record}]\mathbb{P}[\text{obs. }n+1\text{ is }\delta\text{-record}]}\nonumber=\frac{\mathbb{E}[1_{n,\delta}1_{n+1,\delta}]}{\mathbb{E}[1_{n,\delta}]\mathbb{E}[1_{n+1,\delta}]}.
\end{equation}
If the events are independent, then $l_{n}(c,\delta)=1$. Otherwise, values greater or smaller than $1$ indicate positive or negative correlation, respectively. That is, neighbouring $\delta$-records tend to attract or repel each other, if $l_n>1$ or $l_n<1$. 

In order to manipulate $\mathbb{E}[1_{n,\delta}1_{n+1,\delta}]$ we consider the decomposition
\begin{equation}\label{productoindica}
\mathbb{E}[1_{n,\delta}1_{n+1,\delta}]=\mathbb{E}[1_{n,\delta}1_{n+1,\delta}1_{\{Y_n<Y_{n+1}\}}]+\mathbb{E}[1_{n,\delta}1_{n+1,\delta}1_{\{Y_n>Y_{n+1}\}}],
\end{equation}
which, for $\delta<0$, can be written as
\begin{align}\label{eqpnnmas1.1}
\mathbb{E}[1_{n,\delta}1_{n+1,\delta}]&=\int\limits_{-\infty}^{\infty}\!\!\Bigg(  \int\limits_{s-c}^{\infty}\prod_{j=1}^{n-1}F(s+cj-\delta)f(t)dt\! +\int\limits_{s-c+\delta}^{s-c}\prod_{j=2}^nF(t+cj-\delta)f(t)dt \Bigg)\! f(s)ds\nonumber\\
&\!\!\!\!\!\!\!\!=\int\limits_{-\infty}^{\infty}\!\Bigg( ( 1-F(s-c))\prod_{j=1}^{n-1}F(s+cj-\delta)+\!\!\int\limits_{s-c+\delta}^{s-c}\prod_{j=2}^nF(t+cj-\delta)f(t)dt \Bigg)\! f(s)ds,
\end{align}
and, for $\delta\ge0$, 
\begin{align}\label{eqpnnmas1.2}
\mathbb{E}[1_{n,\delta}1_{n+1,\delta}]&=\int\limits_{-\infty}^{\infty}\int\limits_{s-c+\delta}^{\infty}\prod_{j=1}^{n-1}F(s+cj-\delta)f(t)dt  f(s)ds\nonumber\\
&=\int\limits _{-\infty}^{\infty}(1-F(s-c+\delta))\prod_{j=1}^{n-1}F(s+cj-\delta) f(s)ds,
\end{align}
since the second term in \eqref{productoindica} vanishes. 

As for $\mathbb{E}[1_{n,\delta}]$, it is not possible to explicitly compute $\mathbb{E}[1_{n,\delta}1_{n+1,\delta}]$, in general. Nevertheless, it is still possible to describe the behaviour of the dependence index in some particular cases.

\subsection{The Gumbel distribution}
Let $c>0$ and $F$ the Gumbel distribution, as in section \ref{Gumbelsection}. When $\delta<0$ and $n\to \infty$, elementary but lengthy computations yield
\begin{equation*}
\lim_{n\to\infty}\mathbb{E}[1_{n,\delta}1_{n+1,\delta}]=\frac{(e^c-1)^2 (e^c -e^\delta+1)}{(e^c + e^\delta -1)(e^{2c}+ e^\delta -1)}
\end{equation*}
and
\begin{equation*}
l_{\infty}(c,\delta):=\lim_{n\to\infty}l_n(c,\delta)=\frac{(e^c +e^\delta-1) (e^c -e^\delta+1)}{(e^{2c}+ e^\delta -1)}.
\end{equation*}
By differentiating with respect to $c$, we see that $l_\infty(c,\delta)$ is decreasing in $c$ and bounded below by 1, since $\lim_{c\to \infty} l_{\infty}(c,\delta)=1$. With respect to $\delta$ we find that the derivative $\frac{\partial l_{\infty}}{\partial \delta}$ vanishes at
\begin{equation*}
\delta=\log(1-e^{2c}+\sqrt{e^{4c}-e^{2c}}),
\end{equation*}
and then, for any $c$, 
\begin{equation*}
\max_{\delta<0} l_{\infty}(c,\delta)=\frac{2e^{2c}\big( \sqrt{e^{2c}(e^{2c}-1)}-e^{2c}+1 \big)}{\sqrt{e^{2c}(e^{2c}-1)}}=2\big( e^{2c}-\sqrt{2e^{3c}\sinh{(c)}}\big).
\end{equation*}
Note also that $\lim_{\delta\to-\infty}l_{\infty}(c,\delta)=1$.

\begin{figure}
  \centering
  \includegraphics[width=0.5\textwidth]{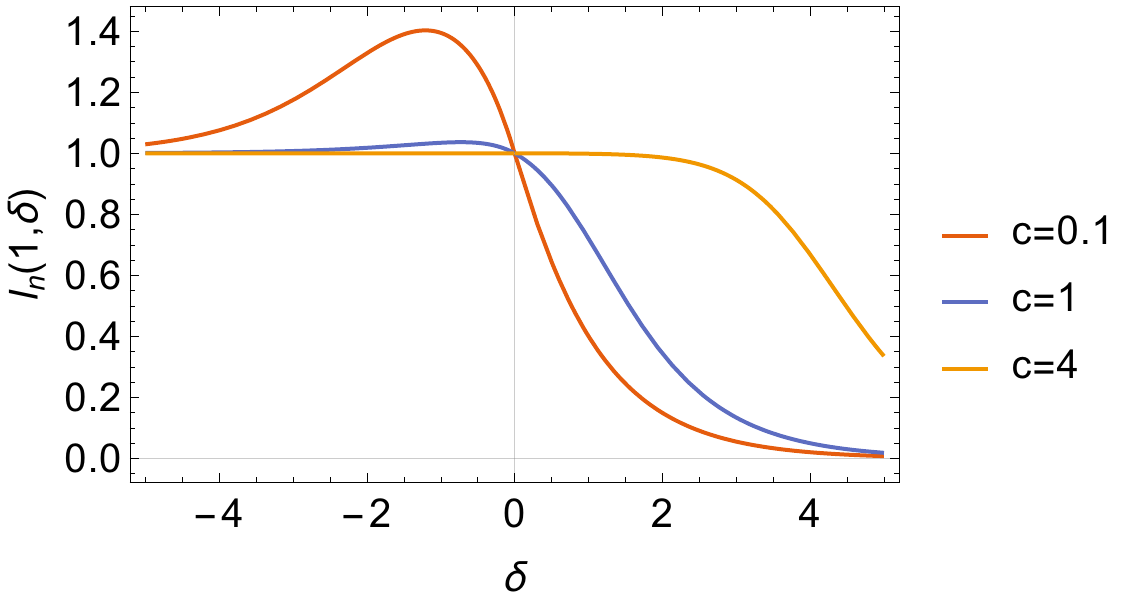}\includegraphics[width=0.5\textwidth]{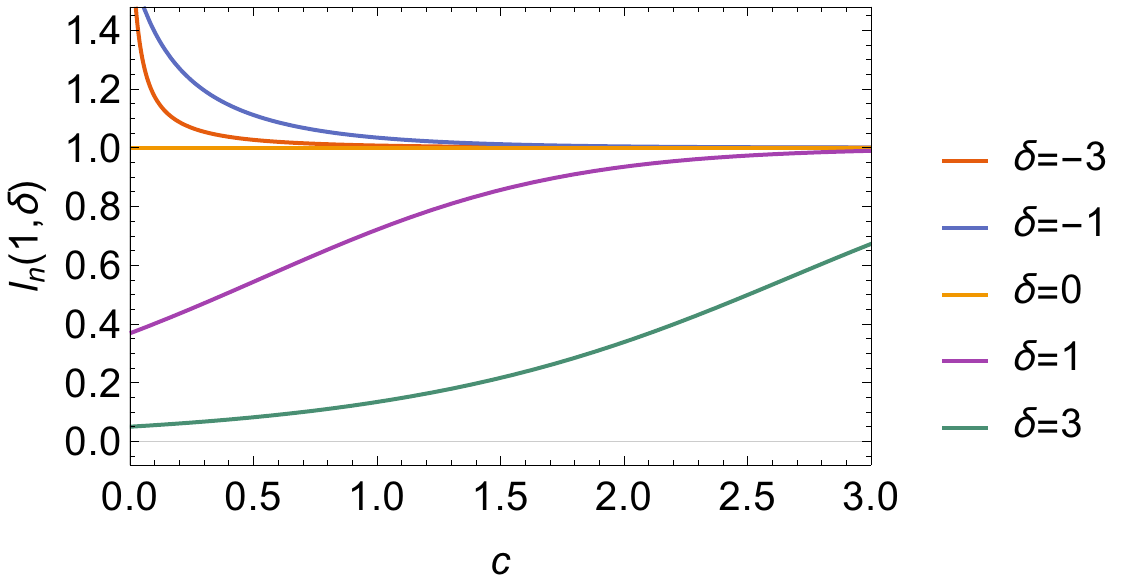}
    \caption{Dependence index $l_\infty(c,\delta)$ for the Gumbel distribution.}
   \label{gumbelPlot}
\end{figure}

For $\delta\ge 0$,
\begin{equation*}
\lim_{n\to \infty}\mathbb{E}[1_{n,\delta}1_{n+1,\delta}]=\frac{e^{c}(e^c-1)^2}{(e^{c}+ e^\delta -1)(e^{c+\delta}-e^c+e^{2c}-e^\delta+e^{2\delta})}
\end{equation*}
and
\begin{equation*}
l_{\infty}(c,\delta)=\frac{e^c(e^c +e^\delta-1) }{e^{c+\delta}-e^c+e^{2c}-e^\delta+e^{2\delta}}.
\end{equation*}
We note that $l_{\infty}(c,\delta)=1$, $\forall c>0$, if $\delta=0$, which results in the asymptotic independence of consecutive record indicators in the LDM. Also, there are no critical points for the index when $\delta\ge 0$. So, in this case  $l_{\infty}(c,\delta)$ is increasing in $c$ with $\lim_{c\to \infty} l_{\infty}(c,\delta)=1$, and  decreasing in $\delta$, with  $\lim_{\delta\to \infty} l_{\infty}(c,\delta)=0$, as can be seen in Figure \ref{gumbelPlot}. Gathering these results, we conclude that $l_{\infty}(c,\delta)>1$ if and only if $\delta<0$. The asymptotic independence for records ($\delta=0$) was proved in \cite{Borovkov1}; 
we have shown here that $\delta$-records attract each other for $\delta<0$ and repel each other for $\delta>0$.

\subsection{The Pareto distribution}
Let $F$ be the Pareto distribution and  $c=1$. The probability of $\delta$-record is given in section \ref{Dagumsection}. Computations of $l_{n}(c,\delta)$ are cumbersome and the explicit expression of $l_n(1,\delta)$ can be found in Appendix \ref{corappendix}. 

We have $\lim_{\delta\to-\infty}l_n(1,\delta)=1$ and $\lim_{\delta\to \infty} l_n(1,\delta)=1-\log(2)\approx 0.3069$, for every $n>1$. Also, $\lim_{n\to \infty} l_n(1,\delta)=\infty$, for all $\delta \in \mathbb{R}$, that is, $\delta$-record-attraction grows  unboundedly, as $n$ increases. Moreover, it can be proved that $ l_n(1,\delta)\sim C\frac{n}{(\log n)^2}$ as $n\to \infty$, where $C$ is a constant depending on $\delta$.

\begin{figure}
  \centering
    \includegraphics[width=0.5\textwidth]{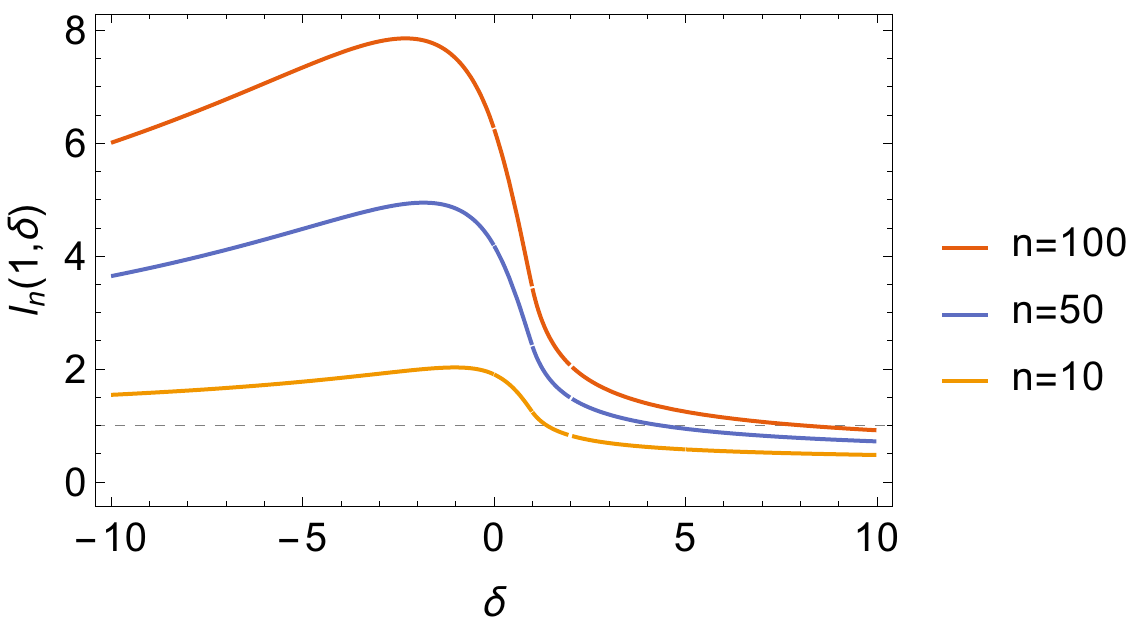}\includegraphics[width=0.5\textwidth]{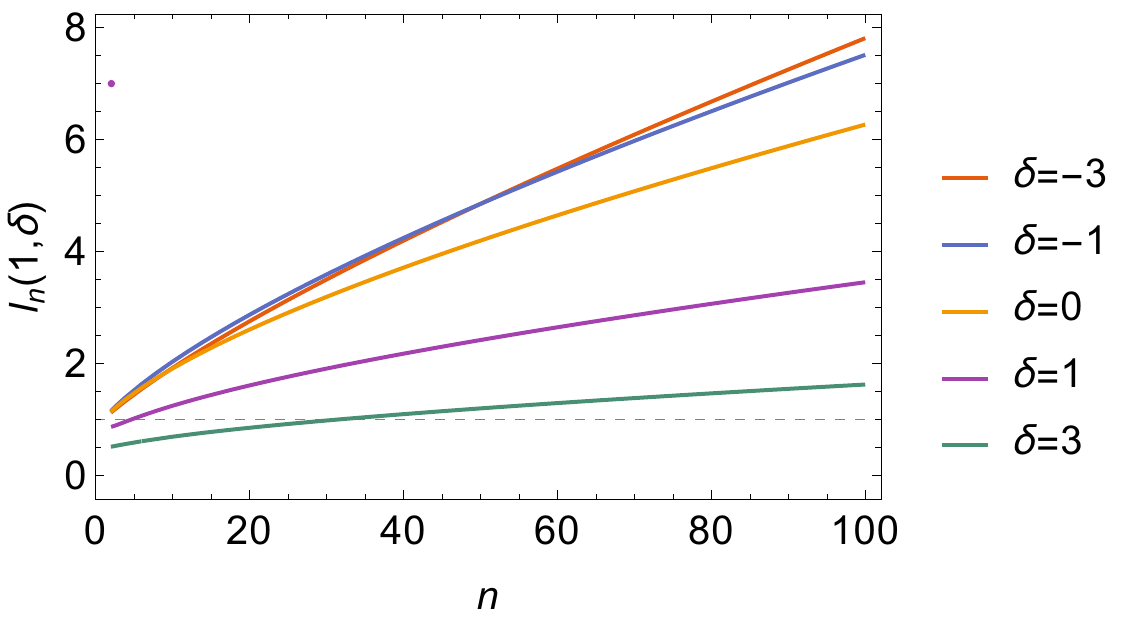}
    \caption{Dependence index $l_n(1,\delta)$ for the Pareto distribution as function of $\delta$ and $n$.}
       \label{paretoPlot}
\end{figure}

The sublinear growth of $l_n (1,\delta)$ as $n$ increases can be observed in the right panel of Figure \ref{paretoPlot}, for different values of $\delta$, as well as the decrease in $\delta$. Also, for fixed $n$ (left panel of Figure \ref{paretoPlot}), there is a negative value of $\delta$ where the correlation  reaches a maximum, as in the Gumbel case. Note that, for negative and small positive values of $\delta$, $l_n (1,\delta)>1$, while, for big values of $\delta$, $l_n(1,\delta)<1$.

\section{Asymptotic behaviour of \texorpdfstring{$N_{n,\delta}$}.}
In sections \ref{seccionanalitica} and \ref{seccion_ejemplos} we have presented properties of the probability that observation $n$ is a $\delta$-record. In this section we analyse the random variable $N_{n,\delta}$, defined as the number of 
$\delta$-records among the first $n$ observations, and study its behaviour as $n\to\infty$. 

Depending on $F$, $c$ and $\delta$, it might be the case that only finitely many  $\delta$-records are observed. We give necessary and sufficient conditions for this to happen. 
On the other hand, if $N_{n,\delta}$ grows to infinity, we investigate if the ratio $N_{n,\delta}/n$ converges (in a certain stochastic sense) to $p_\delta$ and, in that case, how the oscillations of $N_{n,\delta}/n$ around  $p_\delta$ are distributed. 

Recall that, in the classical record model ($c=0$), the number of records $N_{n,0}$ grows to infinity, and there are universal results ensuring that, for any continuous $F$,  $N_{n,0}/\log n$ converges to 1,  almost surely (a.s.) and $(N_{n,0}-\log n )/(\log n)^{1/2}$ has, asymptotically, a standard Gaussian distribution. However, when $\delta\ne0$, results in \cite{Gouet2} and \cite{Gouet22} for the model with $c=0$, show that $N_{n,\delta}$ may grow to a finite limit and, when it diverges, the corresponding limit laws depend both on $\delta$ and $F$. We begin by analyzing the situation where $N_{n,\delta}$ has a finite limit.

\subsection{Finiteness of the total number of $\delta$-records}
Let $N_{\infty,\delta}=\lim\limits_{n\to\infty}N_{n,\delta}$ be the total number of $\delta$-records along the sequence $(Y_n)_{n\ge1}$. In this section we find necessary and sufficient conditions for the finiteness of $N_{\infty,\delta}$ and $\mathbb{E}[N_{\infty,\delta}]$. 

Clearly, these questions are related to the asymptotic behaviour of $p_{n,\delta}$. If $p_\delta>0$, then we can expect $N_{\infty,\delta}=\infty$. On the other hand, if $p_\delta=0$, it may happen that $N_{n,\delta}$ grows sublinearly to $\infty$ or $N_{\infty,\delta}<\infty$. Since, by Theorem \ref{thm:positivity}, the positivity of $p_\delta$ is linked to the finiteness of $\mu^+$, we split the analysis in two cases:
\vspace{4pt}

\noindent {\bf 1.} $\mu^+=\infty$. In this situation, $N_{\infty,\delta}=\infty\ a.s.$ for any $c,\delta\in\mathbb{R}$. 

To check this assertion, we first prove that $M_n:=\max\{Y_1,\ldots,Y_n\}\to\infty$. Observe that $\mu^+=\infty$ implies $x_+=\infty$ and
\begin{equation}
\label{th5}
\displaystyle\sum_{n=1}^{\infty}\mathbb{P}[Y_n>a]= \displaystyle\sum_{n=1}^{\infty}\mathbb{P}[X_n>a-cn] =\sum_{n=1}^{\infty}(1-F(a-cn))=\infty,\ \forall a \in\mathbb{R}.
\end{equation}
From \eqref{th5} and the  second Borel-Cantelli lemma, we conclude that $Y_n>a$ infinitely often (i.o.), for any $a$, and so, $M_n\to\infty$, with probability one. This fact clearly implies $N_{\infty,0}=\infty$. Now, since, for $\delta<0$, $N_{\infty,\delta}\ge N_{\infty,0}$, we get $N_{\infty,\delta}=\infty$. 
On the other hand, for $\delta>0$,  the event
 \begin{equation*}\{X_n+(c-\delta)n>\!\!\max_{1\le j\le n-1}\{X_j+(c-\delta)j\}\}\,\text{ implies }\,
\{X_n+cn>\!\!\max_{1\le j\le n-1}\{X_j+cj\}+\delta\},\end{equation*}
that is, $1_{n,0}(c-\delta)\le1_{n,\delta}(c)$. Therefore,
$N_{\infty,\delta}(c)\ge N_{\infty,0}(c-\delta)=\infty$. 
\vspace{4pt}

\noindent {\bf 2.}  $\mu^+<\infty$. We distinguish three scenarios depending on the sign of $c$.  

If $c>0$, we first assume $x_+-x_->\delta-c$. In this case, we have $p_\delta>0$ and  $N_{\infty,\delta}=\infty$ is an immediate consequence of the law of large numbers in Theorem \ref{slln} below. If $x_+-x_-\le\delta-c$, only the first observation will be a $\delta$-record as shown in section \ref{posp}, so $N_{\infty,\delta}=1$.

If  $c=0$ and $\delta\le 0$, then $N_{\infty,\delta}=\infty$, since $N_{\infty,\delta}\ge N_{\infty,0}=\infty$. If $c=0$ and $\delta>0$, the situation is more complicated. In fact, $N_{\infty,\delta}<\infty$ if and only if 
\begin{equation*}
\int_{0}^{\infty} \frac{1- F(x+\delta)}{(1-F(x))^2}f(x)dx<\infty,
\end{equation*}
which is also equivalent to $\E[N_{\infty,\delta}]<\infty$. This is shown in Proposition \ref{geometricr} of the Appendix, by relating this question to the counting process  of geometric records, as studied in \cite{Gouet22}.

If $c<0$, we proceed as in \eqref{th5} to obtain 
\begin{equation*}
\displaystyle\sum_{n=1}^{\infty}\mathbb{P}[Y_n>a]=\sum_{n=1}^{\infty}\mathbb{P}[X_1>a-cn]<\infty,\quad \forall a \in\mathbb{R},\label{step2thm5}
\end{equation*}
where the last inequality follows from $\mu^+<\infty$. Thus, the first Borel-Cantelli lemma ensures that $\mathbb{P}[Y_n>a\ \text{i.o.}]=0$, for all $a\in\mathbb{R}$, so $Y_n\to-\infty$. Then,  there exists a random variable $N<\infty$ such that $\lim_{n\to\infty} M_n=M_N$ and, consequently,  $N_{\infty,\delta}<\infty$. In this case, we can also prove that $\mathbb{E}[N_{\infty,\delta}]<\infty$; see Proposition \ref{propcnegativo}.

Summarizing the above, we give a complete characterization of the (almost sure) finiteness of the number of $\delta$-records in the next theorem.
\begin{theorem}\label{finiteness}
$N_{\infty,\delta}<\infty\ a.s.$ if and only if one of the following conditions holds
\begin{enumerate}
\item $c<0$ and $\mu^+<\infty$,
\item $c=0$, $\delta>0$ and $\int_{0}^{\infty} \frac{1- F(x+\delta)}{(1-F(x))^2}f(x)dx<\infty$,
\item $c>0$ and $x_+ -x_- \le \delta -c$.
\end{enumerate} 
Moreover, $N_{\infty,\delta}<\infty$ a.s.  if and only if $\mathbb{E}[N_{\infty,\delta}]<\infty$.
\end{theorem}

\begin{remark}Theorem \ref{finiteness} answers a conjecture posed in \cite{Franke2}, stating that the expected number of records ($\delta=0$) in the LDM, with negative trend, remains finite, based on the observed exponential decay of $p_n$, in a particular case. We have shown that the conjecture holds  if  only if $\mu^+<\infty$.
\end{remark}

\subsection{Growth of $N_{n,\delta}$ to infinity.}\label{strongconvergence}
We now turn our attention to the case $N_{\infty,\delta}=\infty$. More precisely, we are interested in the convergence of the proportion of $\delta$-records to $p_\delta$. For records ($\delta=0$) it was shown in \cite{Ballerini-Resnick-1} and \cite{Ballerini-Resnick-2} that $N_{n,0}/n\to p_0$ and that oscillations of $N_{n,0}$ around $p_0$ are asymptotically Gaussian. 

We show here that these results carry over to the case of $\delta\ne0$ but leave the proof for sections \ref{proofslln} and \ref{prooftcl} of the Appendix. As in the aforementioned works, we assume $\mu^+<\infty$ and $c>0$ and, additionally, that $x_+-x_->\delta-c$. Note that, by Theorems \ref{thm:positivity} and \ref{finiteness},  we have $p_\delta>0$ and $N_{\infty,\delta}=\infty$.

\begin{theorem}\label{slln} Assume $\mu^+<\infty$, $c>0$ and $x_+-x_->\delta-c$. Then, as $n\to \infty$, \begin{itemize}
\item[(a)]
$N_{n,\delta}/n\to p_\delta$ a.s. and $\,\mathbb{E}[N_{n,\delta}/n] \to  p_\delta$. 
\item[(b)] If, additionally, $\int_{0}^{\infty}x^2f(x)dx<\infty$, then $\sqrt{n}(N_{n,\delta}/n-p_\delta)\stackrel{\cal D}{\rightarrow} N(0,\sigma_\delta^2)$,  where $\stackrel{\cal D}{\rightarrow}$ stands for convergence in distribution and $\sigma_\delta^2$ is defined in \eqref{varianza}.
\end{itemize}
\end{theorem}

As it can be seen in the proof of Theorem \ref{slln} (a) in the Appendix, the assumption on independence of the  $X_n$ can be relaxed to stationary and ergodic and prove that $N_{n,\delta}/n\to\mathbb{E}[1_{0,\delta}^*]$, defined in \eqref{defestr}. This is useful because it allows to deal with a wider range of scenarios, including stationary ARMA processes. Note, however, that  $\mathbb{E}[1_{0,\delta}^*]$ could differ from $p_\delta$ in \eqref{eq:pdelta}.

\section{Illustration}\label{ilus}

We present a practical application of Theorem \ref{slln} to a real dataset of temperatures, where convergence to the stationary regime is seen for quite small values of $n$. As pointed out in the introduction, the LDM has been used by \cite{Wergen3,Wergen2} to model temperature data in the framework of climate-change.

\begin{figure}
  \centering
    \includegraphics[height=10cm]{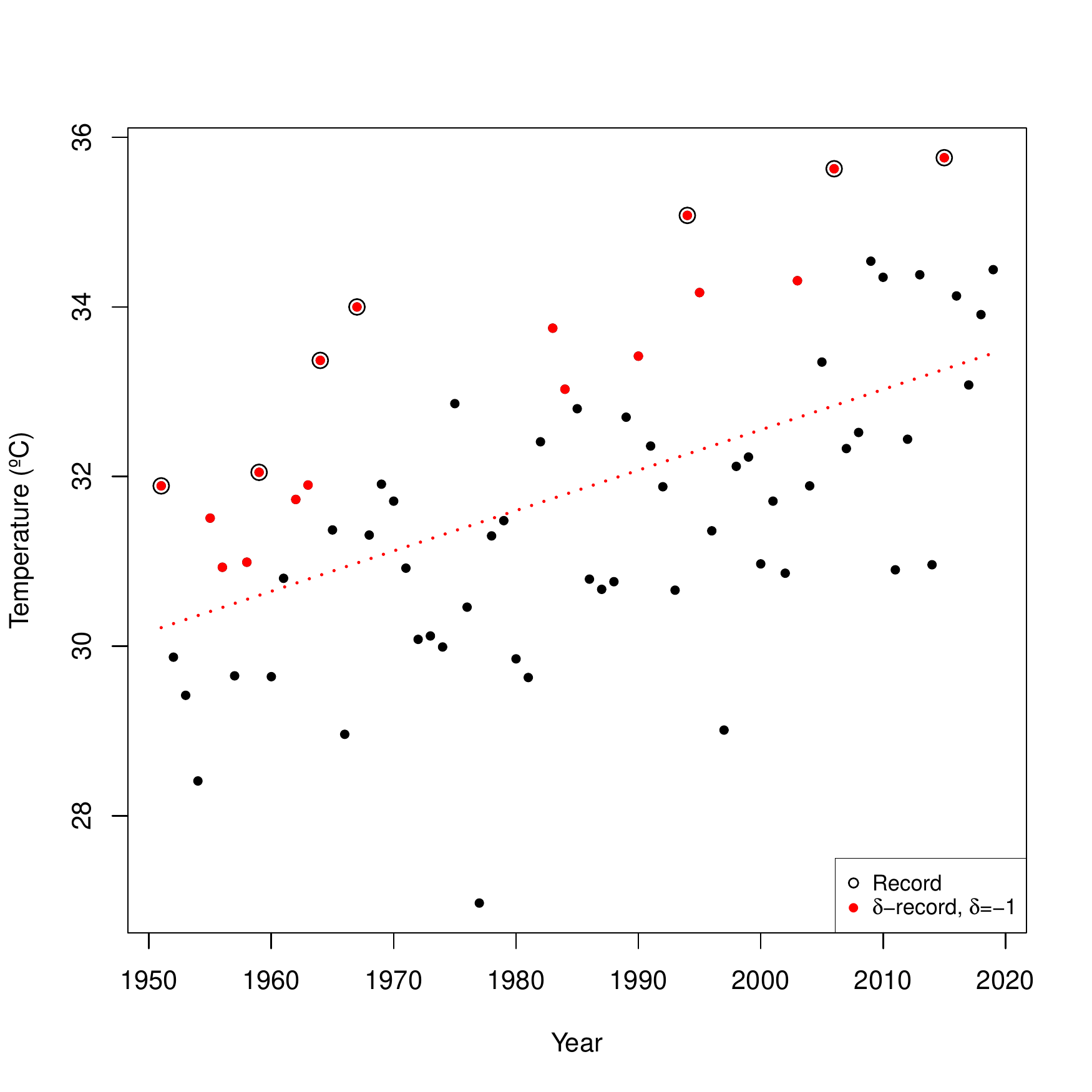}
  \caption{Monthly Mean of Maximum Temperature in July, 1951-2019 in Zaragoza (Spain).}
  \label{fig:plot}
\end{figure}

Our dataset consists of means of daily maximum temperatures (in degrees Celsius), for every month of July, from $1951$ to $2019$, in the city of Zaragoza, Spain. See Figure \ref{fig:plot} for a data plot. The least squares line fitted to the data (in dotted red), reveals a gradual increase of the maximum temperatures over time.
 
 For $\delta$-records we choose the value $\delta=-1$, which is arbitrary and does not respond to any specific reason, other than interpretability of the example. Note that a year will have a $\delta$-record temperature if the maximum average temperature in July is a record or if it is at a distance smaller than $1^\circ C$ from the current maximum. In this framework, we find that $17$ out of the $69$ observations are $\delta$-records (coloured in red), and $7$ of them are records (with circle).
 
The simple linear model for the temperature takes the form
\begin{align}
T_t =\beta_0 +\beta_1t +\varepsilon_{t},\label{model}
\end{align}
where $T_t$ is the temperature of year $t$, and $\varepsilon_{t}$ the error term. The results of the least-square estimators of the coefficients and their $p$-values (assuming Gaussian errors) are shown in Table \ref{ModelResults}.
\begin{table}
	\footnotesize
\begin{center}
  \begin{tabular}{ | c | c | c | c |}
    \hline
      Coefficient & Estimate & Std.Error & p-value \\ \hline
    $\beta_0$ & -62.659 & 18.172 & 0.00098 \\ \hline
    $\beta_1$ & 0.0476 & 0.00915 & 2.04e-06  \\
    \hline
  \end{tabular}
  \caption{Regression analysis estimations for the temperature data.}
  \label{ModelResults}
  \end{center}
\end{table}
In addition, we find an adjusted-$R^2$ of $0.2769$. 

The hypothesis $\beta_1=0$ is clearly rejected, using the Student t-test. Moreover, the estimate of $\beta_1$, which represents the average increment of mean maximum temperatures by year, agrees well with previous estimates of the summer warming trend in Europe, see \cite{Wergen3,Wergen2}.

\begin{figure}
  \centering
    \includegraphics[height=12cm]{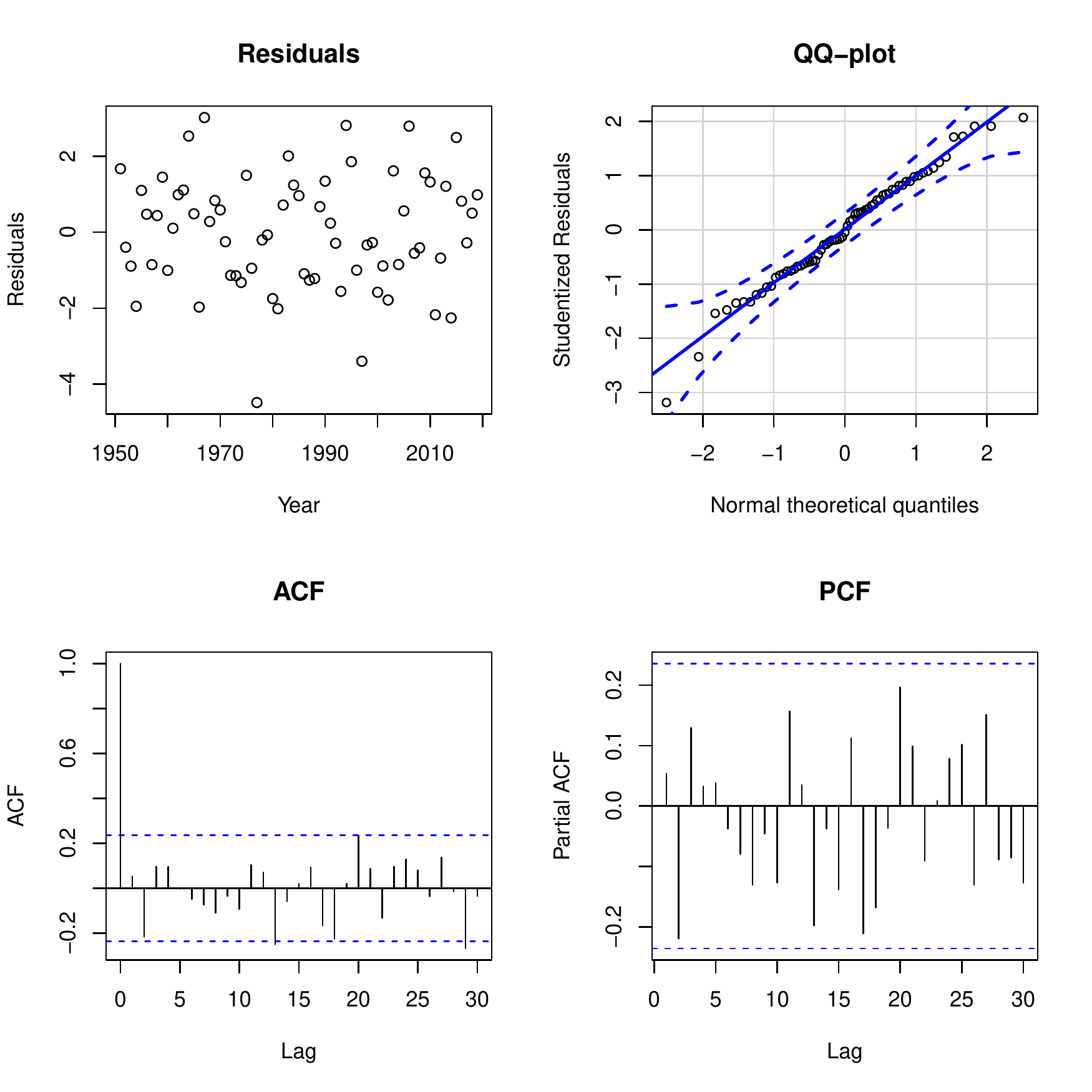}
  \caption{Diagnostic plots of the regression model. Top Left: Residuals vs year. Top Right: Quantile-Quantile of the residual with the normal distribution. Bottom Left: Autocorrelation Function. Bottom Right: Partial Autocorrelation Function.}
  \label{fig:validacion1}
\end{figure}

Figure \ref{fig:validacion1}, along with a $p$-value greater than $0.1$ for the KPSS test of stationarity, and of $0.58$ for the Shapiro-Wilk test applied to the standard regression residuals, indicate that the hypothesis of independent normal residuals is appropriate.

At this point we consider the analogy between the regression model (\ref{model}) and the LDM. First note that the intercept $\beta_0$ is irrelevant when counting $\delta$-records. On the other hand, $\beta_1$ has the role of the trend parameter $c$, whose estimate is $\hat{\beta}_1=0.0476$, as seen in Table \ref{ModelResults}. Finally, the $X_n$ in the LDM are represented by the errors $\varepsilon_t$, which we assume to be zero mean iid.  Note that, for applying Theorem \ref{slln}, there is no need to assume any specific form of the distribution of the $X_n$. 

\begin{figure}
  \centering
    \includegraphics[height=9cm]{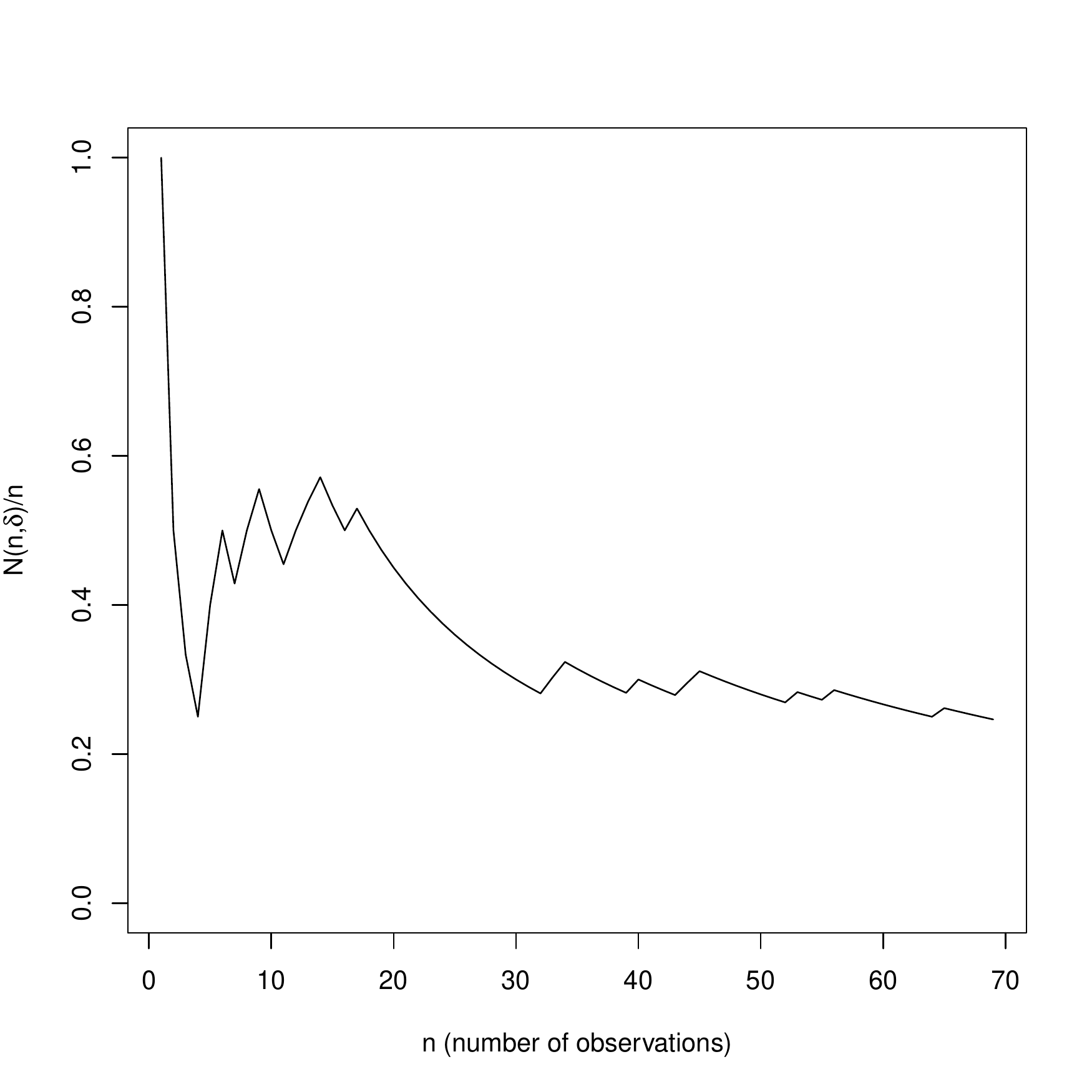}
  \caption{Evolution of the $\delta$-record rate for the temperature data.}
  \label{fig:rate}
\end{figure}

Now, since 17 out of 69 observations were identified as $\delta$-records, it is natural to estimate $p_\delta$ by the empirical record rate, that is,
\begin{equation*}
\hat{p}_\delta=n^{-1}N_{n,\delta}=17/69\approx 0.2464.
\end{equation*} 
Figure \ref{fig:rate} illustrates how the empirical $\delta$-record rate evolves with each extra observation and how it seems to stabilize around a constant value, as  predicted by Theorem \ref{slln}(a).

Concerning the asymptotic normality (Theorem \ref{slln} (b)), we need to estimate the variance $\sigma_\delta^2$, defined in \eqref{varianza}. To that end we propose the estimator 
\begin{equation}
\tilde{\sigma}^2_ \delta=\tilde{\gamma}_{n,\delta}(0)+2\sum_{k=1}^m \tilde{\gamma}_{n,\delta}(k),\label{sigmaest}
\end{equation}
where $m$ is a given natural number and
\begin{equation*}
\tilde{\gamma}_{n,\delta}(k)=n^{-1}\sum_{j=1}^{n-k}(1_{j,\delta}-n^{-1}N_{n,\delta})(1_{j+k,\delta}-n^{-1}N_{n,\delta}),\label{gammaest}
\end{equation*}
The estimator in \eqref{sigmaest} is a version of an estimator proposed in  \cite{Ballerini-Resnick-2}, adapted here to deal with $\delta$-record data. By slightly changing the proof in \cite{Ballerini-Resnick-2}, we can prove convergence of $\tilde{\sigma}^2_\delta$ to $\sigma^2_\delta$, as $n\to\infty$ (consistency), under the condition $m(n)=O(n^{1/2})$.

In order to apply formula \eqref{sigmaest}, we must choose $m$, of order $\sqrt{n}$. In our case, $n=69$ so we take $m=8$, to obtain the estimate $\tilde\sigma_\delta^2=0.337$. Similar values were computed with $m=6,7$. Therefore, from Theorem \ref{slln}(b), $N_{n,\delta}$ is approximately Gaussian, with mean 17 and variance 23.25 ($0.337\times69$). 

\begin{figure}
  \centering
    \includegraphics[height=9cm]{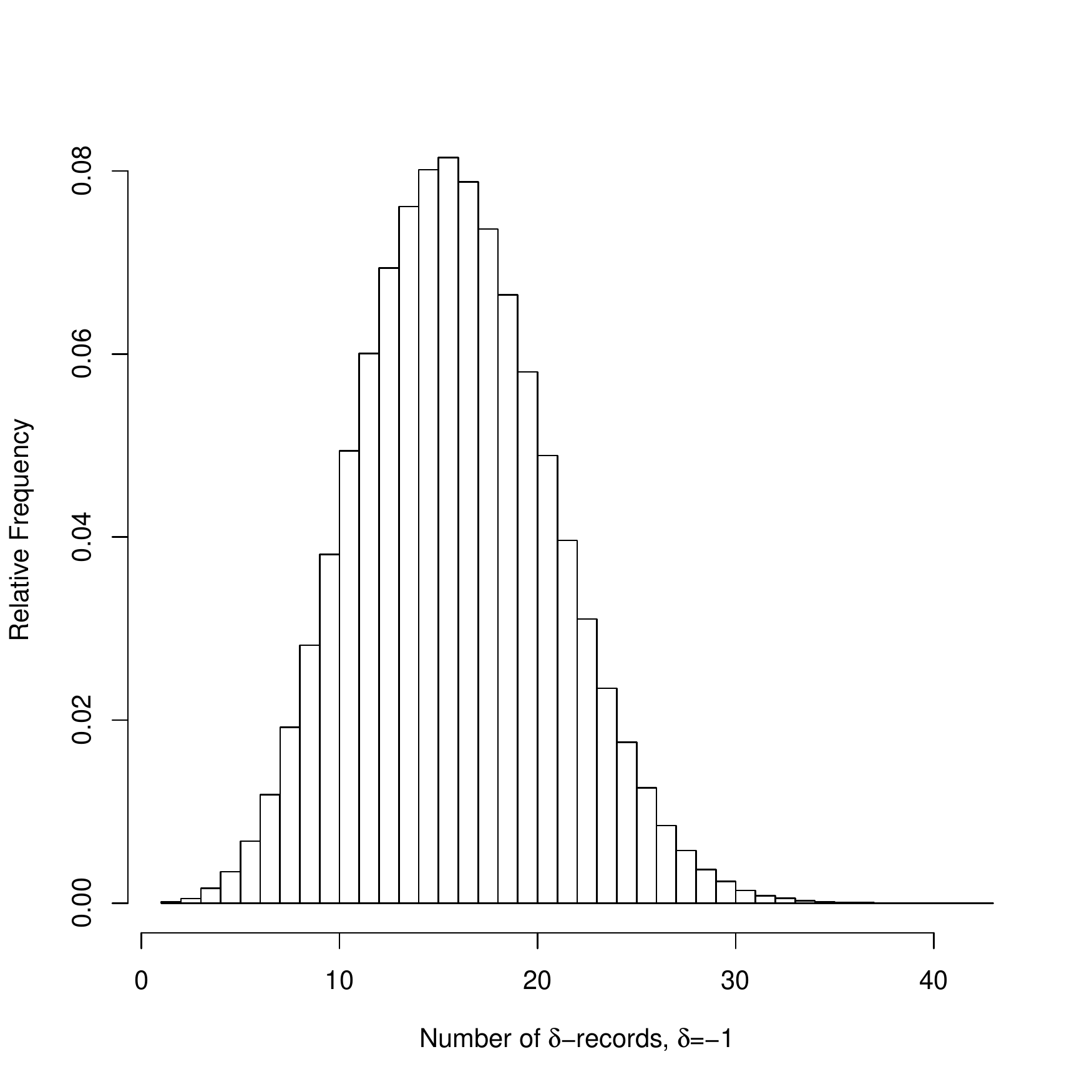}
  \caption{Histogram of the total number of $\delta$-records for the adjusted regression model ($10^6$ iterations of $69$ observations).}
  \label{fig:histogram}
\end{figure}

For assessing the goodness of fit, we simulate the adjusted model \eqref{model} $10^6$ times, and compute the value of $N_{69,\delta}$.  Figure \ref{fig:histogram} summarizes the total number of $\delta$-records obtained at each of the $10^6$ simulations. The histogram has a Gaussian shape, so the convergence in Theorem \ref{slln}(b) to the Gaussian distribution seems to be fast. Moreover, the 0.025 and 0.975 quantiles of the normal distribution $N(17,23.25)$ are, respectively, 7.54 and 26.45. The 0.025 and 0.975 empirical quantiles from the simulated data are 8 and 26, showing an excellent  fit to the theoretical (asymptotic) distribution.

 As a conclusion, we see that empirical results and theory are in very close agreement. This means that, even with a small sample, the approximations in Theorem \ref{slln} are good, at least for the model considered.

\section{Concluding remarks}
In this paper we have studied the behaviour of $\delta$-records in the LDM. We have analysed the asymptotic probability of $\delta$-records,  the dependence between $\delta$-record events and the limiting distribution of the number of $\delta$-records among the first $n$ observations.

The behaviour of the asymptotic probability of $\delta$-records shows similarities with the case of records ($\delta=0$); for instance, for positive $c$, $p_\delta(c)>0$ if and only if $\mu^+<\infty$, regardless the value of $\delta$ (except for the trivial case $\delta\ge x_+-x_-+c$, where no $\delta$-records are observed). We also find that $p_\delta(c)$ is a continuous function of $\delta$ for every $c$, while, as a function of $c$, it is continuous for every $c\ne0$, and a discontinuity arises at $c=0$, if  $x_+<\infty$ and $\delta<0$. This differs from records where $p_0(c)$ is a continuous function of $c$.

We have described in detail the probability of $\delta$-record in two examples. For the Gumbel distribution, an explicit expression for $p_\delta(c)$ is found, showing that it decreases with $\delta$, as a logistic function of $-\delta$. For the cases studied in the Dagum family of distributions, we have $p_\delta(c)=0$, for every $\delta,c$, since $\mu^+=\infty$. 
For this family, we investigate if the speed of convergence of $p_{n,\delta}(c)$ to 0, as $n\to\infty$, depends on $\delta$ or not. Since  random variables $X_n$, with $\mu^+=\infty$, may produce large values, which provoke abrupt changes in record values,  we can expect that $\delta$ values close to 0 have negligible impact and so, $p_{n,\delta}(c)/p_{n,0}(c)\to 1$. This happens in the case $c=0$, where the number of $\delta$-records grows at the same speed as the number of records, when the $X_n$ are heavy-tailed. However, we find that, for some distributions in the Dagum family, $p_{n,\delta}(c)/p_{n,0}(c)\to a\ne1$.

Parameter $\delta$ has a clear impact in the qualitative behaviour of correlations of $\delta$-record events. First, the expression of the limiting correlation is different for $\delta\ge0$ and $\delta<0$. For the Gumbel distribution, where record indicators are independent \cite{Borovkov1}, dependence appears when $\delta\ne0$; in fact, $\delta$-records in this distribution attract each other for $\delta<0$ and repel each other, for $\delta>0$. 
For distributions with power law tails, it is known, for $c>0$, that correlations between records are positive and increase with $n$; see \cite{Franke1}. We have studied the Pareto distribution with $c=1$, and obtained that, while the correlations are positive (and increasing in $n$) for negative, zero and small positive values of $\delta$, they are negative for big values of $\delta$. In fact, for each $n$, the limiting correlation index, as $\delta\to\infty$, is 0.3069.

Another interesting finding of the paper is about the behaviour of the random variable $N_{n,\delta}(c)$. We completely solve the question of finiteness of $N_{\infty,\delta}(c)$, that is, if there is a finite number of $\delta$-records along the infinite sequence of observations. We show that this cannot happen for $c>0$, for any $\delta$ (except if the condition $x_+-x_-<\delta-c$ holds). It cannot happen either when $c<0$ and the underlying random variables $X_n$ have an infinite right-tail mean. This last fact solves a problem posed in \cite{Franke2}, where the authors conjectured that, in the presence of a negative trend, the expected number of records in the whole sequence is finite.

In the case $c>0$ we analyse the asymptotic behaviour of the random variable $N_{n,\delta}$, which grows to infinity. We give a law of large numbers, showing that the ratio $N_{n,\delta}/n$ converges to $p_\delta(c)$ and that its asymptotic distribution is Gaussian, finding the explicit expression of its normalizing constants, which can be estimated from observed data. This result was already known for records and has been applied to different problems, such as athletic records \cite{Ballerini-Resnick-1,Ballerini-Resnick-2} and climate change \cite{Wergen3,Wergen2}. We have illustrated the limiting result for $N_{n,\delta}$ with a set of real data of temperatures of the city of Zaragoza (Spain), showing a good agreement between the theoretical asymptotic results and the observed data in the example. In fact, even for this relatively short series (69 data), the distribution of the number of $\delta$-records is close to the theoretical limiting Gaussian distribution. 

Our results open the door to the use of $\delta$-records for statistical applications in the LDM.  It has been shown that $\delta$-records perform better than records in statistical inference, using trend-free data \cite{Gouet0,Gouet2,Gouet222}, so we expect that their use in the LDM is also advantageous.

\section{Acknowledgements}
This research was funded by project PIA AFB-170001, Fondecyt grant 1161319 and project MTM2017-83812-P of MICINN. The authors are members of the research group Modelos Estoc\'{a}sticos of DGA. M. Lafuente acknowledges the support by the FPU grant, funded by MECD.

We are thankful to professors J. Abaurrea, A.C. Cebri\'an and J. As\'in, from the University of Zaragoza (Spain), for kindly providing us with the temperature data used in this paper.

\section{Appendix}
\subsection{Continuity of $p_\delta(c)$}
\begin{prop}\label{prop:continuidad1} $\prod_{i=1}^\infty F(x+ci-\delta)$, as a function of $c$ 
is continuous at $c\in\mathbb{R}\setminus\{0\}$, for every $x\in(x_-,x_+)$, $x\ne x_-+\delta-c$.
\end{prop}
\textbf{Proof.}
	Let  $(c_n)_{n\ge1}$ be a real sequence converging to $c>0$. We show that
\begin{equation}\label{convint}\prod_{i=1}^\infty F(x+c_ni-\delta)\to\prod_{i=1}^\infty F(x+ci-\delta),\end{equation} as $n\to\infty$, for fixed $x\in(x_-,x_+)$, $x\ne x_-+\delta-c$.

Let $x\in(x_-,x_+)$ be such that $x<x_-+\delta-c$ (this can only happen if $x_->-\infty$ and $\delta-c>0$). In this case $F(x+c-\delta)=0$, so the right-hand side (rhs) of \eqref{convint} is 0. Also, since $c_n\to c$, $F(x+c_n-\delta)=0$, for $n$ large enough,  the left-hand side (lhs) of \eqref{convint} is also 0 and \eqref{convint} is proved.

Let now $x>x_-+\delta-c$, then $F(x+ci-\delta)>0$ for all $i\ge1$. Let $\epsilon>0$ such that $x+c-\epsilon-\delta>x_-$ and let $n_0\ge1$, such that $\vert c_n-c\vert<\epsilon$, for all $n\ge n_0$. We have, for $n\ge n_0$,
\begin{equation*}-\log F(x+c_ni-\delta)\le- \log F(x+(c-\epsilon)i-\delta).\end{equation*}
Since $x>x_-+\delta-(c-\epsilon)$ and $\mu^+<\infty$, we have $-\sum_{i=1}^\infty\log F(x+(c-\epsilon)i-\delta)<\infty$, so the dominated convergence theorem yields
\begin{equation*}\sum_{i=1}^\infty\log F(x+c_ni-\delta)\to\sum_{i=1}^\infty\log F(x+ci-\delta),\end{equation*}
as $n\to\infty$, so \eqref{convint} also holds for $x>x_-+\delta-c$. 
Finally, for $c<0$, we have $\prod_{j=1}^\infty F(x+cj-\delta) = 0, \forall x\in \mathbb{R}$, since $F(x+cj-\delta) \to 0$, as $j\to\infty$.
\hfill$\square$

\subsection{Finiteness of the number of $\delta$-records}
\begin{prop}\label{geometricr}
Let $c=0$ and $\delta>0$. The following conditions are equivalent: 
	\begin{itemize}
		\item[(a)] $N_{\infty,\delta}<\infty$,
		\item[(b)] $\E[N_{\infty,\delta}]<\infty$,
		\item[(c)] \begin{equation*}\label{condicioninfini}
\int_{0}^{\infty} \frac{1- F(x+\delta)}{(1- F(x))^2}f(x)dx<\infty.
\end{equation*}	
	\end{itemize}
\end{prop}
\textbf{Proof.} It is clear that $Y_n$ is a $\delta$-record if and only if $e^{X_n}>e^\delta\max\{e^{X_1},\dots,e^{X_{n-1}}\}.$
That is, if the $n$-th observation in the sequence $(e^{X_n})_{n\ge1}$ is a geometric record, with parameter $k=e^\delta$, according to \cite{Gouet22}. In section 2.1.1 of that paper, it is shown that the total number of geometric records, in a sequence of iid random variables, with cdf $G$, is finite if and only if 
\begin{equation}\label{condiciong}
\int_1^\infty\frac{ 1- G(kx)}{(1-G(x))^2}dG(x)<\infty.\end{equation}
Moreover, in section 2.3.4 of that paper, it is shown that \eqref{condiciong} is equivalent to the finiteness of the expectation of the total number of geometric records.
Since $G(x)=F(\log(x))$, the result is proved.
\hfill$\square$

 In the rest of the Appendix, we use the operator $\bigvee$ to denote the maximum. Then, for instance, $\bigvee_{i=1}^nY_i=\max\{Y_1,\ldots,Y_n\}$.
\begin{lemma}\label{lemanegativo}
\begin{enumerate}
\item If $c<0$, $x_->-\infty$ and $\mu^+<\infty$, then $\mathbb{E}[N_{\infty,\delta}]<\infty,\ \forall \delta\in\mathbb{R}.$
\item Let $\tilde{X}_1$ be a random variable with cdf $G$, and $(\tilde{X}_n)_{n\ge2}$ an iid sequence, independent of $\tilde{X}_1$, with common cdf $F$, such that $G(x)\le F(x),\ \forall x$. Let $\tilde{Y}_n=\tilde{X}_n+cn, n\ge1$. Then, if $c<0$,  $\mathbb{E}[\sum_{j=1}^\infty 1_{\{\tilde{Y}_j>\vee_{i=1}^{j-1}\tilde{Y_i}+\delta\}}]\le\mathbb{E}[N_{\infty,\delta}]$.
\end{enumerate}

\end{lemma}
\textbf{Proof.} $(i)$  First we bound $p_{n,\delta}(c)$ as follows
\begin{align*}p_{n,\delta}(c)&=\int_{-\infty}^\infty \prod_{j=1}^{n-1}F(x+cj-\delta)f(x)dx\nonumber\\
&=\int_{-\infty}^\infty \prod_{j=1}^{n-1}F(x+cj-\delta)1_{\{x+c(n-1)-\delta>x_-\}}f(x)dx\nonumber\\
&=\int_{x_--c(n-1)+\delta }^\infty \prod_{j=1}^{n-1} F(x+cj-\delta)f(x)dx\\ 
&\le 1-F(x_--c(n-1)+\delta).
\end{align*}
So, $\sum_{j=1}^n p_{n,\delta}\le \sum_{j=1}^n (1-F(x_--c(j-1)+\delta))$ yielding
\begin{equation*}
\mathbb{E}[N_{\infty,\delta}]\le \sum_{j=1}^\infty (1-F(x_--c(j-1)+\delta))<\infty,
\end{equation*}
since $\mu^+<\infty.$

$(ii)$ It suffices to check that the $\delta$-record probability for the $\tilde{Y}_n$ fulfills
\begin{align*}
\mathbb{E}[1_{\{\tilde{Y}_j>\vee_{i=1}^{j-1}\tilde{Y_i}+\delta\}}]&=\int_{-\infty}^\infty G(x+c-\delta) \prod_{i=2}^{j-1}F(x+ci-\delta)f(x)dx\nonumber\\
&\le \int_{-\infty}^\infty \prod_{i=1}^{j-1}F(x+ci-\delta)f(x)dx =p_{j,\delta}(c). \quad \square
\end{align*}

\begin{prop}\label{propcnegativo}
If $c<0$ and $\mu^+<\infty$, then $\mathbb{E}[N_{\infty,\delta}]<\infty$.
\end{prop}
\textbf{Proof.}
It suffices to consider $\delta<0$, since the number of $\delta$-records is decreasing with $\delta$. Also, we take $x_-=-\infty$ as, otherwise, the result follows from Lemma \ref{lemanegativo} $(i)$. Moreover, since there exists $c_1 \in \mathbb{R}$ such that $\mathbb{P}(X_n+c_1>0)>0$, and  the number of $\delta$-records is the same for the sequences $Y_n=X_n+cn$ and $\tilde{Y}_n=X_n+cn+c_1$, we assume without loss of generality that $\mathbb{P}(X_n>-\delta)>0.$

Let $N=\inf\{n\in\mathbb{N}\mid X_n>-\delta\}$, then  $N$ is a geometric random variable and
\begin{equation*}
N_{\infty,\delta}=\sum_{j=1}^N 1_{j,\delta}+\sum_{j=N+1}^\infty 1_{j,\delta}
=\sum_{j=1}^N 1_{j,\delta}+\sum_{j=N+1}^\infty 1_{j,\delta}1_{\{X_{j}>0\}}.
\end{equation*}
For $j>N$, let $\tilde{1}_{j,\delta}=1_{\{X_j>\vee_{i=N}^{j-1}(X_i+c(i-j)+\delta)\}}1_{\{X_j>0\}}$, then
\begin{equation*}
1_{j,\delta}1_{\{X_j>0\}}=1_{\{X_j>\vee_{i=1}^{j-1}(X_i+c(i-j)+\delta)\}}1_{\{X_j>0\}}\le\tilde{1}_{j,\delta}.
\end{equation*}
Note that the $\tilde{1}_{j,\delta}$, defined for $j>N$, are the $\delta$-record indicators of the sequence $\{X_N,X_{N+1}1_{\{X_{N+1}>0\}}+c,X_{N+2}1_{\{X_{N+2}>0\}}+2c,\dots\}$. Now, taking expectations we have
\begin{align*}
\mathbb{E}[N_{\infty,\delta}]\le \frac{1}{\mathbb{P}(X_1>\-\delta)}+\sum_{i=1}^\infty\mathbb{E}[\tilde{1}_{i,\delta}]<\infty,
\end{align*}
since the last sum is bounded by Lemma \ref{lemanegativo} $(ii)$. \hfill$\square$
 
\subsection{Proof of law of large numbers for \texorpdfstring{$N_{n,\delta}(c)$}{TEXT} } \label{proofslln}

Define the bilateral LDM, as in \eqref{lineardrift}, but letting $n\in\mathbb{Z}$ instead of $n\in\mathbb{N}$.
Associated to this model, define, for $n\in\mathbb{Z}$,
\begin{equation}\label{defestr}
M_n^*=\max\{Y_i:\ i\le n\}, \qquad 1_{n,\delta}^*=1_{\{Y_n>M_{n-1}^*+\delta\}},
\end{equation}
and, for $n\in\mathbb{N}$,
\begin{equation*}
N_{n,\delta}^*=\sum_{k=1}^n 1^*_{k,\delta}.
\end{equation*}

\noindent\textbf{Theorem} \ref{slln} (a).
Let $c>0$, $\mu^+<\infty$. Then $N_{n,\delta}(c)/n\to p_\delta(c)$ a.s.  as $n\to \infty$.

\noindent\textbf{Proof.}
It is clear that
\begin{equation*}
\lim_{n\to\infty}\mathbb{P}[Y_n>a]=\lim_{n\to\infty}\mathbb{P}[X_n>a-cn]= 1,\ \forall a\in \mathbb{R}, \label{step1thm5} 
\end{equation*}
thus  $Y_n\to\infty\  $ and $M_n \to \infty$  a.s. Also, since $\mu^+<\infty$, it is known by a Borel-Cantelli argument that $M_0^*< \infty$ a.s. Gathering these facts, we know that $\exists \ 0<N<\infty$ a.s. such that $ 1_{N,0}^*=1$ almost surely. From the definition of $1^*_{n,0}$, given $n\in\mathbb{N}$ we  have $1_{n,0}\ge1^*_{n,0}$, and so $1_{N,0}=1$ a.s., entailing  $M_n^*=M_n$ and $1_{n,\delta}=1^*_{n,\delta}\ a.s.\ \forall n>N$. So,
\begin{align*}\label{sumigual}
\sum_{k=N+1}^\infty 1_{k,\delta}=\sum_{k=N+1}^\infty 1^*_{k,\delta}\ a.s.
\end{align*}
Also, we know that $1_{n,\delta}^*$ is a strictly stationary and ergodic sequence. Applying Birkhoff's Theorem we have
\begin{equation*}
\frac{N_{n,\delta}^*}{n}=\frac{1}{n}\sum_{k=1}^n 1^*_{k,\delta} \to \mathbb{E}[1_{0,\delta}^*]\ a.s.
\end{equation*}
Now,  let $(a_n)_{n\ge1}$ be a real sequence diverging to $\infty$. Then
\begin{equation*}
\left\lvert \frac{N_{n,\delta}-N_{n,\delta}^*}{a_n}\right\lvert \le \left\lvert  \frac{N}{a_n}\right\lvert \to 0\ a.s.
\end{equation*}
since $N$ does not depend on $n$. Finally, since $\left\lvert \frac{N_{n,\delta}-N_{n,\delta}^*}{n}\right\lvert \to 0\ a.s.$ and $\frac{N^*_{n,\delta}}{n} \to \mathbb{E}[1_{0,\delta}^*]\ a.s.$, we have $\frac{N_{n,\delta}}{n} \to \mathbb{E}[1_{0,\delta}^*]\ a.s.$ Finally, $\mathbb{E}[1_{0,\delta}^*]$ can be written as the rhs in (\ref{eq:pdelta}), yielding $\mathbb{E}[1_{0,\delta}^*]=p_\delta(c)$.

\subsection{Proof of central limit theorem for \texorpdfstring{$N_{n,\delta}(c)$}. }\label{prooftcl}

A proof of Gaussian convergence for the number of $\delta$-records, based on the ideas in \cite{Ballerini-Resnick-1}, is not straightforward. The main problem arises when considering the joint probability of two observations being $\delta$-records. While in the case of records this quantity can be explicitly written as follows
\begin{align*}
\mathbb{E}[1_{i,0}1_{i+m,0}]=\int_{-\infty}^{\infty}\prod_{k=1}^{i-1}F(y+ck)\int_{y-cm}^{\infty}\prod_{j=1}^{m-1}F(s+cj)f(s)dsf(y)dy,
\end{align*}
in the setting  $\delta\ne 0$ there is no such analytical expression. In order to solve this problem we introduce the following general bounds, which do not depend on the specification of the model for the sequence $(Y_n)_{n\ge1}$.

\begin{prop} \label{bounds}
Let $(Y_k)_{k\in \mathbb{Z}}$ be a sequence of random variables and consider the events $A=\Big\{\bigvee\limits_{k=-\infty}^{i-1} Y_{k} + \delta < Y_i\Big\}$, B=$\Big\{ \bigvee\limits_{k=i+1}^{i+m-1} Y_{k} + \delta < Y_{i+m}\Big\}$,  $C=\{Y_i - \delta < Y_{i+m}\}$ and $E=\{Y_i + \delta < Y_{i+m}\}$. Then, if $\delta\le0$,
\vspace{4pt}

\noindent a1)  $\mathbb{P}[A\cap B\cap C]\leq\mathbb{E}[1_{i,\delta}^*1_{i+m,\delta}^*]$ and 
\vspace{4pt}

\noindent a2) $\mathbb{P}[A\cap B\cap E]\geq\mathbb{E}[1_{i,\delta}^*1_{i+m,\delta}^*]$.
\vspace{4pt}

Also, if $\delta\ge 0$,

\noindent b) $\mathbb{P}[A\cap B\cap E]=\mathbb{E}[1_{i,\delta}^*1_{i+m,\delta}^*]$.
\end{prop}
\textbf{Proof}.
\textit{a1)} Note that $1_{j,\delta}^*$ is the indicator of  $D_j=\Big\{\bigvee\limits_{k=-\infty}^{j-1} Y_{k} + \delta < Y_j\Big\}, j=i,i+m$. Then we must show that $A\cap B\cap C\subseteq D_i\cap D_{i+m}$.

First, it is clear that $A= D_i$. Also, observe that $C\subseteq E$ and that $A\cap C\subseteq \Big\{\bigvee\limits_{k=-\infty}^{i-1} Y_{k} + \delta < Y_{i+m}\Big\}$, since $\delta\le0$. From the inclusions above we have
$$ A\cap B\cap C\subseteq \Big\{\bigvee\limits_{k=-\infty}^{i-1} Y_{k} + \delta < Y_{i+m}\Big\}\cap E \cap B=D_{i+m}$$ and the conclusion follows.

\textit{a2)} Trivial.\\

\textit{b)} It is clear that $D_i\cap D_{i+m}\subseteq A\cap B\cap E$ and that $A\cap B\cap E\subseteq D_i$, because $A=D_i$. Also, since $\delta\ge0$, we have $A\cap E\subseteq \Big\{\bigvee\limits_{k=-\infty}^{i-1} Y_{k} + \delta < Y_{i+m}\Big\}$, so 
$$ A\cap B\cap E\subseteq \Big\{\bigvee\limits_{k=-\infty}^{i-1} Y_{k} + \delta < Y_{i+m}\Big\}\cap E \cap B=D_{i+m},$$ which completes the proof. \hfill$\square$

Note that, although it is unnecessary in our setting, the reverse \textit{a1)} inequality also holds for $\delta\ge0$. Under the assumptions of the LDM, the lhs of the first two bounds in the previous proposition have analytical expressions. The strategy to prove Gaussian convergence  is  to work with the corresponding bounds of $\mathbb{E}[1_{i,\delta}1_{i+m,\delta}]$, which are shown to be tight enough to achieve our purpose.
So, with this result we slightly modify the necessary bounds and rebuild the martingale approach in \cite{Ballerini-Resnick-1}, to prove convergence to the Gaussian distribution.
\vspace{4pt}

\textbf{Theorem} \ref{slln} (b).\label{thm: clt}
Suppose that $\int_0^\infty x^2f(x)dx< \infty $ and let $c>0,\ \delta \in \mathbb{R}$, such that $p_\delta>0$. Then, as $n\to \infty$,
\begin{center}
$\sqrt{n}(n^{-1}N_{n,\delta}-p_\delta(c))\stackrel{\cal D}{\rightarrow} N(0,\sigma_\delta^2(c)),$
\end{center}
where \begin{equation}\label{varianza}\sigma_\delta^2 = p_\delta-p_\delta^2+2\sum_{m=1}^{\infty}(\mathbb{E}[1_{i,\delta}^* 1_{i+m,\delta}^*]-p_\delta).\end{equation}

\textbf{Proof.} 
For simplicity, we only consider the case $\delta\le 0$ since the case $\delta > 0$ is analogous. We assume $-2\delta<x_+$ as, otherwise, we can define $X_n'=X_n+(-3\delta-x_+)$, $n\ge1$; the number of $\delta$-records in both models is the same and $-2\delta<x_+'$, where $x_+'$ is the right-end point of $X_n'$.

The proof is split into several steps.

\noindent\textbf{1)} 
We claim that
\begin{align}
0\le p_{n,\delta}-p_{\delta}\le c^{-1}\int_{c(n-1)/2-\delta}^\infty(1-F(s))ds+F(-\delta)^{\floor{(n-1)/2}}.\label{cltpdbound}
\end{align}
The first inequality follows from
\begin{align*}
p_{n,\delta}-p_{\delta}= \int_{-\infty}^\infty\left(\prod\limits_{j=1}^{n-1}F(y+cj-\delta)-\prod\limits_{j=1}^{\infty}F(y+cj-\delta)\right)f(y)dy\ge 0.
\end{align*}
For the second, let $u=\prod\limits_{j=1}^{n-1}F(y+cj-\delta)$ and $v=\prod\limits_{j=1}^{\infty}F(y+cj-\delta)$. Then, from the elementary inequality $u-v\le u-uv$, we have 
\begin{equation}
\label{eq:integral}
p_{n,\delta}-p_{\delta}\le\int_{-\infty}^{\infty}u(1-v)f(y)dy.
\end{equation}
The integral in the rhs of \eqref{eq:integral} is split into two terms $A, B$, that we bound. Let $A=\int_{-\infty}^{-c(n-1)/2}u(1-v)f(y)dy$ and $B=\int_{-c(n-1)/2}^{\infty}u(1-v)f(y)dy$, then
\begin{align}
\label{eq:A}
A&\le \int_{-\infty}^{-c(n-1)/2}\prod\limits_{j=1}^{n-1}F(-c(n-1)/2+cj-\delta)f(y)dy\nonumber\\
&\le\prod\limits_{j=1}^{n-1}F(c(j-(n-1)/2)-\delta)\nonumber\\
&\le \prod\limits_{j=1}^{\floor{(n-1)/2}}F(c(j-(n-1)/2)-\delta)\nonumber\\
&\le \prod\limits_{j=1}^{\floor{(n-1)/2}}F(-\delta)=F(-\delta)^{\floor{(n-1)/2}}.
\end{align}
For $B$ we have
\begin{align}
\label{eq:B}
B&\le  \int_{-c(n-1)/2}^{\infty}\left(1-\prod\limits_{j=n}^{\infty}F(y+cj-\delta)\right)f(y)dy\nonumber\\
&\le \int_{-c(n-1)/2}^{\infty}\sum\limits_{j=n}^{\infty}(1-F(y+cj-\delta))f(y)dy\nonumber\\
&\le \int_{-c(n-1)/2}^{\infty}\left( \int_{z=n-1}^{\infty}(1-F(y+cz-\delta))dz \right) f(y)dy\nonumber\\
&\le \int_{-c(n-1)/2}^{\infty}\left( c^{-1} \int_{-c(n-1)/2+c(n-1)-\delta}^{\infty}(1-F(s))ds \right) f(y)dy\nonumber\\
&\le c^{-1}\int_{c(n-1)/2-\delta}^\infty(1-F(s))ds.
\end{align} 
So, from \eqref{eq:A} and \eqref{eq:B}, \eqref{cltpdbound} holds.

\noindent\textbf{2)} Let $r_{m,\delta}=\mathbb{E}[1_{i,\delta}^*1_{i+m,\delta}^*]$, which is well defined since it does not depend on $i$. We bound $r_{m,\delta}$ by applying Proposition \ref{bounds} as follows:
\begin{align*}
r_{m,\delta}&=\mathbb{P}\left[Y_i,Y_{i+m}\mbox{ are }\delta\mbox{-records}\right]\nonumber\\
&=\mathbb{P}\left[Y_i>\bigvee_{l<i}Y_l+\delta,Y_{i+m}>\bigvee_{l<i+m}Y_l+\delta\right]\nonumber\\
&\le \mathbb{P}\left[Y_i>\bigvee_{l<i}Y_l+\delta,Y_{i+m}>\bigvee_{l=1}^{m-1}Y_{i+l}+\delta,Y_{i+m}>Y_i+\delta\right]\nonumber\\
&=\iint_{y<s+cm-\delta}\prod\limits_{j=1}^{\infty}F(y+cj-\delta)\prod\limits_{i=1}^{m-1}F(s+ci-\delta)f(s)dsf(y)dy.
\end{align*}

If $r_{m,\delta}\ge p_{\delta}^2$, we apply the Fubini-Tonelli theorem, as well as the triangle inequality, to obtain
\begin{align*}
|r_{m,\delta}- p_{\delta}^2|\le& \left|\iint_{y<s+cm-\delta}\prod\limits_{j=1}^{\infty}F(y+cj-\delta)\prod\limits_{i=1}^{m-1}F(s+ci-\delta)f(s)dsf(y)dy-p_\delta^2\right|\nonumber\\
\le& A+B,
\end{align*}
where $$A=\int_{-\infty}^\infty\prod\limits_{j=1}^{\infty}F(y+cj-\delta)\left|\int_{-\infty}^{\infty}\prod\limits_{i=1}^{m-1}F(s+ci-\delta)f(s)ds-p_\delta\right|f(y)dy\nonumber$$ and 
$$B=\int_{-\infty}^\infty\prod\limits_{j=1}^{\infty}F(y+cj-\delta)\int_{-\infty}^{y-cm+\delta}\prod\limits_{i=1}^{m-1}F(s+ci-\delta)f(s)dsf(y)dy.$$

Since variables are separated in $A$ and applying the first step of this proof
\begin{align}\label{pcltA}
A&\le \int_{-\infty}^\infty\prod\limits_{j=1}^{m-1}F(s+cj-\delta)f(s)ds-p_\delta\nonumber\\
&\le c^{-1}\int_{c(n-1)/2-\delta}^\infty(1-F(s))ds+F(-2\delta))^{\floor{(m-1)/2}}.
\end{align}
While for $B$ we have
\begin{align}\label{pcltB}
B=& \int_{-\infty}^{cm/2}\prod\limits_{j=1}^{\infty}F(y+cj-\delta)\int_{-\infty}^{y-cm+\delta}\prod\limits_{i=1}^{m-1}F(s+ci-\delta)f(s)dsf(y)dy\nonumber\\
&+ \int_{cm/2}^\infty\prod\limits_{j=1}^{\infty}F(y+cj-\delta)\int_{-\infty}^{y-cm+\delta}\prod\limits_{i=1}^{m-1}F(s+ci-\delta)f(s)dsf(y)dy\nonumber\\
\le& \int_{-\infty}^{-cm/2+\delta}\prod\limits_{j=1}^{m-1}F(s+cj-\delta)f(s)ds
\int_{-\infty}^{cm/2+}\prod\limits_{j=1}^{\infty}F(y+cj-\delta)f(y)dy\nonumber\\
&+\int_{cm/2}^\infty \prod\limits_{j=1}^{\infty}F(y+cj-\delta)f(y)dy\nonumber\\
\le& \prod\limits_{j=1}^{m-1}F(-cm/2+cj)+1-F(cm/2)\nonumber\\
\le& F(-2\delta)^{\floor{(m-1)/2}}+1-F(cm/2).
\end{align}
Analogously, applying the corresponding bound in Proposition \ref{bounds}, we arrive at the same conclusion if $r_{m,\delta}\le p_\delta$ via (\ref{pcltA}) and (\ref{pcltB}), so
\begin{equation}
|r_{m,\delta}-p_\delta|\le c^{-1}\int_{c(m-1/2-\delta)}^\infty(1-F(s))ds+2F(-2\delta)^{\floor{(m-1)/2}}+1-F(cm/2)\label{rmbound}.
\end{equation}

\noindent\textbf{3)} Since $\int_0^\infty x^2f(x)dx<\infty$, it is easy to check, from (\ref{rmbound}), that the series $\sum_{m=1}^\infty |r_{m,\delta}-p_\delta|$ converges; for $F(-2\delta)^{\floor{(m-1)/2}}$ convergence holds since $F(-2\delta)<1$.

\noindent\textbf{4)} Using the strategy in the proof of Theorem \ref{slln} (a), we get the  following convergence in distribution
\begin{equation}
\sqrt{n}(n^{-1}N_{n,\delta}-n^{-1}N_{n,\delta}^*)\stackrel{\cal D}{\rightarrow}0.\label{normtransfer}
\end{equation}

\noindent\textbf{5)} Theorem 5.2 in \cite{Hall} is applied to the $N_{n,\delta}^*$ in order to transfer the asymptotic normality to $N_{n,\delta}$,  as a consequence of  (\ref{normtransfer}). This martingale result guarantees convergence to the Gaussian distribution, if the next two conditions hold: 
\begin{enumerate}
\item$ \sum_{k=1}^\infty \mathbb{E}[\xi_{k,\delta}\mathbb{E}[\xi_{l,\delta}|\mathcal{M}_0]]$ converges $\forall l\ge 0.$
\item$\lim_{l\to\infty} \sum_{k=K}^\infty \mathbb{E}[\xi_{k,\delta}\mathbb{E}[\xi_{l,\delta}|\mathcal{M}_0]]=0$ uniformly in  $K\ge1.$
\end{enumerate}
where $\xi_{k,\delta}=1_{k,\delta}^*-p_\delta$ and $\mathcal{M}_0$ is a certain sub-$\sigma$-algebra of events of the original probability space (see \cite{Hall}, page 128, for details). Moreover we have

$$\lim_{n\to\infty} n^{-1}\mathbb{E}\left[\left(\sum_{i=1}^n \xi_{i,\delta}\right)^2\right]=\sigma_\delta^2.$$

Given that the hypothesis $\delta\ne 0$ does not imply any extra difficulty in the application of this theorem, we omit the verification of these two conditions since the adaptation of this part of the proof is forthright following the lines of \cite{Ballerini-Resnick-1}.
\hfill$\square$

\subsection{Correlations in the Pareto Distribution}\label{corappendix}
The $\delta$-record probability is given in \eqref{probpareto}.
For $\mathbb{E}[1_{n,\delta}1_{n+1,\delta}]$ and $n>2$, we use (\ref{eqpnnmas1.1}) for $\delta<0$ and (\ref{eqpnnmas1.2}) for $\delta\ge0.$

\noindent{\bf 1.} Let $a=n-\delta$, $A=(\delta -2) (\delta  (1-a)+(n-1) \log a) (n \log (a+1)-\delta a),$  
$$B=-\big(\delta ^3 (n-2)+\delta -2 n^3-2 \delta ^2 (n^2-2)+\delta  (n-1) (n+5)n+n+1\big) \log (a+1),$$
 $$C=(a-1) \log (a+1- \delta)-(\delta -2)a\left(\delta  (a-1)^2-(n-1)a \log (4a)\right)+(1-a) \log ((a-\delta +1 ) (a+1)).$$ Then, if $\delta<0$, 
 \begin{equation*}
 l_n(1,\delta)=\frac{B+C}{A}.
 \end{equation*}

   \noindent{\bf 2.} Let $a=n-\delta$, $A=(\delta -1)^2 (\delta -a) (\delta  (1-a)+(n-1) \log a) (-\delta  a+n \log (a+1))$,
   $$B=(a- \delta ) \left((\delta -1) (\delta ^2 (a-1)+(\delta -1) (n-1)\log \left(\frac{a- \delta+1}{(2-\delta)a}\right)\right),$$
   $$C=-\log (2-\delta ) (\delta 
   (\delta +2)-2 \delta  n+n-1)\big)+(\delta -1)^2 (n-1) \log (a- \delta +1).$$
    Then, if $0<\delta<1$,
      \begin{equation*}
l_{n}(1,\delta)=\frac{a^2 (B+C)}{A}.
\end{equation*}

\noindent{\bf 3.}  If $\delta=1$, 
   \begin{equation*}
    l_{n}(1,1)=\frac{(n-1)^2 ((n-2) n-2 (n-1) \log (n-1))}{2 (n-2) (-n+(n-1) \log (n-1)+2) (-n+n \log (n)+1)}.
   \end{equation*}
    \noindent{\bf 4.} Finally, let $a=n-\delta$, $A_1=(\delta +\log \delta-n \log \delta-n+(n-1) \log (n-1)+1)(\delta -1)^2 (\delta -a)$, $$A_2= (\delta -n \log \delta   -n+n \log n),$$
   $$B=(\log \delta ) \big(2 \delta  (\delta ^2+2 \delta -1)+(2 \delta -1) n^2-5 \delta ^2 n+n\big),$$
   $$C=(\delta -1)^2 (n-1) \log(n-1)-(a-1) \big((\delta -1) (\delta -a)+(2 \delta -1) \log (2 \delta -1) (a-1)\big).$$
   Then,  if $\delta>1$, and $\delta \notin \{n/2,n,n+1\}$ (otherwise, the values of the index are given by the continuous extension at these points),
   \begin{equation*}
   l_n(1,\delta)=\frac{a^2(B+C)}{A_1A_2}.
   \end{equation*}

\end{document}